\documentclass[12pt]{amsart}
\usepackage{latexsym, amssymb, amsmath, amscd}
 \usepackage{eucal}

    \newtheorem{rema}{Remark}[section]
    \newtheorem{propo}[rema]{Proposition}

   \newtheorem{theo}[rema]{Theorem}
   \newtheorem{def-theo}[rema]{Definition-Theorem}
 \newtheorem{conj}[rema]{Conjecture}
   \newtheorem{defi}[rema]{Definition}
    \newtheorem{lemma}[rema]{Lemma}
    \newtheorem{corol}[rema]{Corollary}
     \newtheorem{exam}[rema]{Example}
  \newtheorem{rmk}[rema]{Remark}

	\newcommand{\nno}{\nonumber}

	\newcommand{\p}{\partial}

 \newcommand{\res}{\mbox{\rm Res}}

 \newcommand{\pf}{{\it Proof:}\hspace{2ex}}
 
 \newcommand{\epfv}{\hspace{1em}$\Box$\vspace{1em}}

\newcommand{\bC}{{\mathbb C}}
\newcommand{\bZ}{{\mathbb Z}}
\newcommand{\bR}{{\mathbb R}}
\newcommand{\bQ}{{\mathbb Q}}
\newcommand{\bN}{{\mathbb N}}

\newcommand{\cA}{{\mathcal A}}

\newcommand{\bD}{{\mathbb D}}
\newcommand{\cU}{{\mathcal U}}

\newcommand{\cM}{{\mathcal M}}
\newcommand{\cN}{{\mathcal N}}

\newcommand{\cB}{{\mathcal B}}

\newcommand{\cC}{{\mathcal C}}

\newcommand{\cz}{\bC[z]}
\newcommand{\czz}{\bC[z^{-1}, z]}

\newcommand{\Kz}{K[z]}
\newcommand{\cAz}{\cA[z]}
\newcommand{\cAzz}{\cA[z^{-1}, z]}

\newcommand{\im}{\mbox{Im\,}}

\newcommand{\IC}{{\bf IC }}

\newcommand{\BQ}{\begin{eqnarray}}
\newcommand{\EQ}{\end{eqnarray}}
\newcommand{\BQn}{\begin{eqnarray*}}
\newcommand{\EQn}{\end{eqnarray*}}

\title[The Image Conjecture and the Mathieu Conjecture]
{Generalizations of the Image Conjecture and the Mathieu Conjecture}

  \author{Wenhua Zhao}      
    \date{\today}

\begin{document}

\begin{abstract}
We first propose a generalization of the image conjecture 
\cite{IC} for the commuting differential operators related with classical orthogonal polynomials. We then show that the non-trivial case of this generalized image conjecture is equivalent to a variation of the Mathieu conjecture \cite{Ma} from integrals of $G$-finite functions over reductive Lie groups $G$ to integrals of polynomials over open subsets of $\bR^n$ with any positive measures. Via this equivalence, the generalized image conjecture can also be viewed as a natural variation of Duistermaat and van der Kallen's theorem \cite{DK} on Laurent polynomials with no constant terms. To put all the conjectures above in a common setting, 
we introduce what we call the {\it Mathieu} subspaces of associative algebras. We also discuss some examples of Mathieu subspaces from other sources and derive some general results on this newly-introduced notion.
\end{abstract}

   

 
\thanks{The author has been partially supported 
by NSA Grant R1-07-0053}

\bibliographystyle{alpha}
    \maketitle


\renewcommand{\theequation}{\thesection.\arabic{equation}}
\renewcommand{\therema}{\thesection.\arabic{rema}}
\setcounter{equation}{0}
\setcounter{rema}{0}
\setcounter{section}{0}

\section{\bf Introduction}

\subsection{Background and Motivation}

The main motivations and contents of this paper are as follows.

First, in \cite{IC} a so-called {\it image} 
conjecture ({\bf IC}) on images of 
commuting differential operators of polynomial algebras of 
order one with constant leading coefficients
has been proposed. It has also been shown there 
that the well-known {\it Jacobian} conjecture proposed by O. H. Keller \cite{Ke} (See also \cite{BCW} and \cite{E}) and, more generally, the {\it vanishing} conjecture 
\cite{HNP}, \cite{GVC} on differential operators 
(of any order) with constant coefficients,  
are actually equivalent to some special 
cases of the {\bf IC}.

Second, as pointed out in \cite{GVC}, all 
classical orthogonal polynomials in one or 
more variables can be obtained from some 
commuting differential operators of order 
one with constant leading 
coefficients. Unfortunately, most of these 
differential operators are not differential 
operators of polynomial algebras. 
Instead, they are differential  
operators of some localizations of 
polynomial algebras such as Laurent 
polynomial algebras, etc.

Note that, due to their applications in many 
different areas of mathematics such as in ODE, PDE, 
the eigenfunction problems and representation theory, 
orthogonal polynomials have been under intense 
study by mathematicians in the last two centuries. 
For example, in \cite{SHW} published in $1940$, about $2000$ 
published articles on orthogonal polynomials mostly 
in one variable had been included. 
Therefore it will also be interesting to consider the 
\IC for the commuting differential operators 
related with classical orthogonal polynomials.

Unfortunately, the straightforward generalization of 
the \IC from polynomial algebras to their localizations 
does not hold in general. In this paper, we propose another 
generalization of the \IC (See Conjecture \ref{GIC}) 
for the commuting differential operators related with classical orthogonal polynomials. 

We will also show that, under certain conditions, 
the new generalization is actually 
equivalent to a conjecture (See Conjecture \ref{GMC}) on integrals of polynomials over open subsets 
of $B\subset \bR^n$ with any (positive) measures. 
The latter conjecture turns out to be a natural variation of the Mathieu conjecture \cite{Ma} (See Conjecture \ref{MC}) from $G$-finite functions on reductive Lie groups $G$ to polynomial functions over the open subsets 
$B\subset\bR^n$ above. It also can be viewed as a natural variation of  Duistermaat and van der Kallen's theorem \cite{DK} (See Theorem \ref{ThmDK}) on Laurent polynomials with no constant terms. 

To be more precise, let us first introduce the following notion 
which will provide a common ground for all the results and 
conjectures to be discussed in this paper.
 
\begin{defi}\label{M-ideal}
Let $R$ be any commutative ring and $\cA$ a commutative $R$-algebra.
We say that a $R$-subspace 
$\cM$ of $\cA$ is a {\it Mathieu subspace} of $\cA$  
if the following property holds: 
for any $a, b\in \cA$ with $a^m\in \cM$ for any 
$m\ge 1$, we have $a^m b \in \cM$ when $m>>0$, 
i.e. there exists $N\ge 1$ $($depending on $a$ and $b$$)$ 
such that $a^m b \in \cM$ for any $m\ge N$.
\end{defi}   

Note that, any ideal of $\cA$ is automatically 
a Mathieu subspace of $\cA$. But conversely, 
not all Mathieu subspaces are ideals. 
Actually, many Mathieu subspaces are not even closed 
under the product of the ambient algebra $\cA$.
So the new notion can be viewed as a generalization 
of the notion of ideals. For more examples and general results 
on Mathieu subspaces, see Section \ref{S5}. 

The notion is named after Olivier Mathieu due to 
his following conjecture proposed in \cite{Ma}, $1995$. 

\begin{conj} \label{MC} $(${\bf The Mathieu Conjecture}$)$
Let $G$ be a compact connected real Lie group with the 
Haar measure $\sigma$. Let $f$ a complex-valued $G$-finite 
function over $G$ such that $\int_G f^m \, d\sigma=0$ 
for any $m\ge 1$. Then, for any $G$-finite function $g$ 
over $G$, $\int_G f^m g \,  d \sigma =0$ when $m>>0$.
\end{conj}

Note that, in terms of the newly introduced notion of 
Mathieu subspaces, the Mathieu conjecture just claims that 
the $\bC$-subspace of complex-valued $G$-finite functions 
$f$ with $\int_G f d\sigma=0$ is a Mathieu subspace of the 
$\bC$-algebra $\cA$ of complex-valued $G$-finite 
functions over $G$. 

One of the motivations of the Mathieu conjecture is its 
connection with the Jacobian conjecture (See \cite{BCW} and \cite{E}). Actually, Mathieu also showed in \cite{Ma} that his conjecture implies the Jacobian conjecture.  

For later purposes, here we also point out that 
J. Duistermaat and W. van der Kallen \cite{DK} in $1998$ had proved the Mathieu conjecture for the case of tori, which now can be re-stated as follows.

\begin{theo}\label{ThmDK} $(${\bf Duistermaat and 
van der Kallen}$)$ 
 Let $z=(z_1, z_2,$ $..., z_n)$ be $n$ commutative free   
variables and $\cM$ the subspace of 
the Laurent polynomial algebra $\czz$ 
consisting of the Laurent polynomials 
with no constant terms. Then $\cM$ is a Mathieu 
subspace of $\czz$. 
\end{theo}

Another main motivation behind the new notion of Mathieu subspaces is the following so-called {\it image} conjecture ({\bf IC}) proposed recently by the author in \cite{IC} on the images of commuting differential operators of polynomial algebras of order one with constant leading coefficients.

Let $z=(z_1, z_2, ..., z_n)$ be $n$ commutative free  
variables and $\cz$ the algebra of polynomials in $z$ 
over $\bC$. For any $1\le i\le n$, set 
$\p_i\!:=\p/\p z_i$. We say a differential 
operator $\Phi$ of $\cz$ is {\it of order one 
with constant leading coefficients} if 
$\Phi=h(z)+\sum_{i=1}^nc_i\p_i$ for some 
$h(z)\in \cz$ and $c_i\in \bC$. We denote by 
$\bD[z]$ the subspace of all differential 
operators of order one with constant 
leading coefficients. For any subset 
$\cC=\{\Phi_i\,|\, i\in I\}$ of differential 
operators of $\cz$, we set 
$\im\cC\!:=\sum_{i\in I}(\Phi_i \cz)$ 
and call it the {\it image} of $\cC$. 
We say $\cC$ is {\it commuting} if, for any $i, j\in I$, 
$\Phi_i$ and $\Phi_j$ commute with each other.

With the notation fixed above, the \IC can be re-stated as follows.
 
\begin{conj} \label{IC} $(${\bf The Image Conjecture}$)$
For any commuting subset $\cC\subset\bD[z]$, $\im \cC$ 
is a Mathieu subspace of $\cz$.
\end{conj}

Note that the {\bf IC}, the Mathieu conjecture and 
also Conjectures \ref{GIC}--\ref{GMC} mentioned at the beginning of this subsection are all problems on whether or not certain subspaces are Mathieu subspaces. It is also the case for the Jacobian conjecture and, more generally, the vanishing conjecture \cite{HNP}, \cite{GVC} on differential operators (of any order) with constant coefficients via their connections with the \IC (See \cite{IC}). 
Furthermore, we can also include the well-known 
Dixmier conjecture \cite{D} in the list 
since it has been shown, first by Y. Tsuchimoto 
\cite{Ts} in $2005$ and later by A. Belov and M. Kontsevich
\cite{BK} and P. K. Adjamagbo and 
A. van den Essen \cite{AE} in $2007$, that the Dixmier conjecture 
is actually equivalent to the Jacobian conjecture.
The implication of the Jacobian conjecture 
from the Dixmier conjecture was actually proved 
much earlier by V. Kac (unpublished 
but see \cite{BCW}) in $1982$. 

Therefore, it is interesting and important to study Mathieu subspaces separately in a general and abstract setting. 
So we will also discuss more examples of Mathieu subspaces 
from other sources and derive some general 
results on this newly introduced notion 
(See Section \ref{S5}). 

\subsection{Arrangement} 
In Subsection \ref{S4.1}, we first recall some classical 
orthogonal polynomials and their related commuting differential operators (See Examples \ref{COP} and \ref{MJP}). We also fix some notations and summarize some facts that will be needed for the rest of this paper. 

In Subsection \ref{S4.2}, we consider the straightforward generalization of the \IC for the commuting differential operators related with some multi-variable Jacobi orthogonal polynomials but without the constraints on the parameters required by the Jacobi polynomials. We will show in Proposition \ref{P4.2.1} that the straightforward generalization of the \IC does not hold for these differential operators. But, if we generalize the \IC in a different way, we will have a positive answer for these differential operators 
under the constraints on the parameters 
required by the Jacobi polynomials 
(See Corollary \ref{C4.2.3}).

Another purpose of this subsection 
is to explain in a concrete setting 
the main ideas behind the generalization of 
the \IC that will be formulated and discussed 
in Section \ref{S4b}. Some of the results of this subsection 
will also be needed later in Subsection \ref{S4.4}. 

In Subsection \ref{S4.3}, we first formulate a 
generalization (See Conjecture \ref{GIC}) of the {\bf IC} for the differential operators related with orthogonal polynomials, and also a conjecture (See Conjecture \ref{GMC})  
on integrals of polynomials over open subsets 
$B\subset \bR^n$ with any positive measures.  
We show in Proposition \ref{P4.2.3} that the non-trivial case of 
Conjecture \ref{GIC} is actually 
equivalent to some special cases of 
Conjecture \ref{GMC}. We also point 
out that Conjecture \ref{GMC} in some sense 
can be viewed as a natural variation of 
the Mathieu conjecture (See Conjecture \ref{MC}) 
and Duistermaat and van der Kallen's theorem 
(See Theorem \ref{ThmDK}).

In Subsection \ref{S4.4}, we prove some cases of 
Conjectures \ref{GIC} and \ref{GMC}. We also discuss 
a connection of Conjecture \ref{GMC} with the 
{\it polynomial moment problem} which was first proposed 
by M. Briskin, J.-P. Francoise and Y. Yomdin in the series paper 
\cite{BFY1}-\cite{BFY5} and recently was solved by  
F. Pakovich and M. Muzychuk \cite{PM}. 

In Section \ref{S5}, we discuss Mathieu subspaces in the most general setting. Some examples of Mathieu subspaces from other sources will be given and some general results on this newly introduced notion will also be derived. \\

{\bf Acknowledgment}
The author greatly thanks Harm Derksen, Jean-Philippe Furter, Jeffrey C. Lagarias, Leonid Makar-Limanov, Lucy Moser-Jauslin for communications and suggestions on a connection of Conjecture 
\ref{GMC} with the polynomial moment problem. The author is also very grateful to Mitya Boyarchenko for sending the author a sketch of his brilliant but unpublished proof of Theorem \ref{Folk} and to Arno van den Essen for pointing out some errors and misprints of an earlier version of this paper. 
At last but certainly not the least, great thanks also go to Fedor Pakovich for communications on his joint work \cite{PM} with Mikhail Muzychuk.

\renewcommand{\theequation}{\thesection.\arabic{equation}}
\renewcommand{\therema}{\thesection.\arabic{rema}}
\setcounter{equation}{0}
\setcounter{rema}{0}

\section{\bf Differential Operators Related with Classical Orthogonal Polynomials}\label{S4a}

In this section, we first recall in Subsection \ref{S4.1} 
some classical orthogonal polynomials in one or more variables 
and their related differential operators. 
We also summarize in Lemma \ref{L4.1.5} some facts 
that will be needed in later sections. 

The classical reference for one-variable orthogonal polynomials 
is \cite{Sz} (see also \cite{AS}, \cite{C}, \cite{Si1}). 
For multi-variable orthogonal polynomials, 
see \cite{DX}, \cite{Ko} and references therein. 
But here we will essentially follow the presentations 
given in \cite{GVC} in terms of differential operators
of certain localizations of polynomial algebras,  
and emphasize that the related differential operators 
are all commuting differential operators 
of order one with constant leading coefficients.
The presentation in \cite{GVC} for 
multi-variable orthogonal polynomials 
will also be simplified here.

In Subsection \ref{S4.2}, we consider the 
straightforward generalization of the image conjecture 
({\bf IC}) for  some commuting differential operators 
of the Laurent polynomials. Up to changes  
of variables, these differential operators 
are related with some multi-variable
Jacobi orthogonal polynomials (See Example \ref{MJP}). 
We show in Proposition \ref{P4.2.1} that the \IC does 
not hold for these differential operators. But 
it does hold for the same differential operators 
if we generalize the \IC in a different way 
(See Corollary \ref{C4.2.3}).

\subsection{Differential Operators Related with Classical Orthogonal Polynomials}\label{S4.1}

First, let us recall the definition of 
classical orthogonal polynomials. 
In order to be consistent with the traditional notations 
of orthogonal polynomials, in this subsection we will 
 use $x=(x_1, x_2, \dots, x_n)$ instead of 
$z=(z_1, z_2, \dots, z_n)$ to denote free 
commutative variables.

\begin{defi}\label{Def-OP}
Let $B$ be a non-empty open subset of $\bR^n$ and $w(x)$ 
a real valued function defined over $B$ such that 
$w(x)\geq 0$ for any $x\in B$ and 
$0<\int_B w(x)dx <\infty$. 
Assume further that $\int_B f(x) w(x)dx$ 
is finite for any  $f(x)\in \bC[x]$.
A sequence of polynomials 
$\{u_\alpha (x) \,|\, \alpha\in \bN^n \}$ is 
said to be {\it orthogonal} over $B$ 
if 
\begin{enumerate}
\item[$(a)$] $\deg u_\alpha=|\alpha|:=\sum_{i=1}^n k_i$ 
 for any $\alpha=(k_1, k_2, ..., k_n)\in \bN^n$. 
\item [$(b)$] the sequence 
$\{u_\alpha(x)\,|\, \alpha\in \bN^n\}$ forms an orthogonal basis of $\bC[x]$ with respect to the Hermitian form defined by 
\begin{align} \label{H-Form}
(f, g)\!:=\int_B  f \bar g  w(x) \, dx
\end{align}
for any $f, g\in \bC[x]$, where $\bar g$ denotes the complex conjugate of the polynomial $g\in \bC[x]$.
\end{enumerate}
\end{defi}

The function $w(x)$ is called the {\it weight} function.
For all classical orthogonal polynomials, $w(x)$ is smooth over $B$ but might have some singular points over the boundary of $B$ (See Examples \ref{COP} and \ref{MJP} below). When the open set $B\subset \bR^n$ and $w(x)$ are 
clear in the context, we simply call the polynomials 
$u_\alpha(x)$ $(\alpha\in \bN^n)$ in the definition 
above {\it orthogonal polynomials}.
If the orthogonal polynomials 
$u_\alpha(x)$ $(\alpha\in \bN^n)$ 
also satisfy $\int_B |u_\alpha|^2 w(x)dx=1$ 
for any $\alpha\in \bN^n$, we call $u_\alpha(x)$ 
$(\alpha\in \bN^n)$ {\it orthonormal polynomials}.

Note that, if $u_\alpha(x)$ $(\alpha\in \bN^n)$ are orthogonal polynomials, say, as in Definition \ref{Def-OP}, then, 
for any $c_\alpha \in \bC^{\times}$ $(\alpha\in \bN^n)$, 
$c_\alpha u_\alpha$ $(\alpha\in \bN^n)$ 
are also orthogonal polynomials over $B$ with the same weight function $w(x)$. 

An obvious way to construct orthogonal polynomials is to apply 
the Gram-Schmidt process. But, surprisingly, most of classical 
orthogonal polynomials can also be obtained by the following  
so-called Rodrigues' formulas which, 
in terms of the notation as in Definition \ref{Def-OP}, 
can be stated as follows. \\

{\bf Rodrigues' Formula}: {\it There exist some non-zero constants $c_\alpha \in \bR$ $(\alpha\in \bN^n)$ and an $n$-tuple $g(x)=(g_1(x), g_2(x), ..., g_n(x))$ of 
polynomials in $x$ such that}
\begin{align}\label{Rodrigues}
u_\alpha(x)=c_\alpha w(x)^{-1} \frac {\p^{|\alpha|}}{dx^\alpha} 
(w(x)g^\alpha(x)). \\ \nno
\end{align}

Note that, not all orthogonal polynomials 
defined in Definition \ref{Def-OP} can be 
obtained by Rodrigues' formulas.  
For example, the weight function $w(x)$ 
of some orthogonal polynomials may 
not even be differentiable.

Now, for any $1\le i\le n$, let 
\begin{align}\label{L4OP}
\Lambda_i \!:=w(x)^{-1} \left (\frac {d}{dx_i}\right) w(x)
=\frac {d}{dx_i}+w(x)^{-1}\frac {dw(x)}{dx_i}.
\end{align}  
and set   
$\Lambda\!:=(\Lambda_1, \Lambda_2, ..., \Lambda_n)$.
Then, 
by Rodrigues' formula above, we see that the 
orthogonal polynomials 
$\{u_\alpha (x)\,|\,\alpha\in \bN^n \}$ have the form
\begin{align}\label{LaGa}
u_\alpha(x)=c_\alpha \Lambda^\alpha (g^\alpha(x))
\end{align}
for any $\alpha\in \bN^n$.

Note also that the differential operator $\Lambda_i$ in 
Eq.\,(\ref{L4OP}) is a differential operator of order one 
with constant leading coefficients. Furthermore, in the 
multi-variable case, the differential operators 
$\Lambda_i$ $(1\le i\le n)$ commute with one another 
since they are the conjugations of the commuting 
differential operators $\p_i$ $(1\le i\le n)$ 
by the multiplication operator by $w^{-1}(z)$.

Let us look at the following classical orthogonal polynomials.

\begin{exam} \label{COP} \mbox{}

\begin{enumerate}
\item {\bf Hermite Polynomials:} 

$(a)$ $B=\bR$ and the weight function $w(x)=e^{-x^2}$. 

$(b)$ the differential operator $\Lambda$ and the polynomial 
$g(x)$:
\begin{align}\label{HmtEqs}
\begin{cases}
\Lambda  &=\frac d{dx}-2x, \\
g(x)&=1,
\end{cases}
\end{align}

$(c)$ the Hermite polynomials in terms of $\Lambda$ and $g(x)$:
\begin{align*}
H_m(x)=(-1)^m \, \Lambda^m (g^m(x)).
\end{align*}

\item {\bf Laguerre Polynomials:} 

$(a)$ $B=\bR^+$ and 
$w(x)=x^\alpha e^{-x}$ $(\alpha>-1)$.

$(b)$ the differential operator $\Lambda$ and the polynomial $g(x)$:
\begin{align}\label{LgrEqs}
\begin{cases}
\Lambda  &=\frac d{dx}+ (\alpha x^{-1}-1), \\
g(x)&=x,
\end{cases}
\end{align}

$(c)$ the Laguerre polynomials in terms of $\Lambda$ and $g(x)$:
\begin{align*}
L_m(x)=\frac 1{m!}\,  \Lambda^m (g^m(x)).
\end{align*}

\item 
{\bf Jacobi Polynomials:} 

$(a)$ $B=(-1, 1)$ and 
$w(x)=(1-x)^\alpha (1+x)^\beta$ with $\alpha, \beta >-1$.

$(b)$ the differential operator $\Lambda$ and the polynomial $g(x)$:
\begin{align}\label{JP-Lambda}
\begin{cases}
\Lambda &=\frac{d}{dx}-\alpha(1-x)^{-1}+\beta(1+x)^{-1},
 \\
g(x)&=1-x^2.
\end{cases}
\end{align}

$(c)$ the Jacobi polynomials in terms of $\Lambda$ and $g(x)$:
\begin{align}\label{JP-LP}
P^{\alpha, \beta}_m (x)=\frac {(-1)^m}{2^m m!}\, \Lambda^m g^m(x).
\end{align}

\item {\bf Classical Orthogonal Polynomials over Unit Balls:}

$(a)$ $B={\mathbb B}^n=\{x\in \bR^n\,|\, ||x||<1\}$ 
and the weight function
\begin{align*}
w_\mu (x)=(1-||x||^2)^{\mu - 1/2},
\end{align*}
where $||\cdot||$ denotes the usual Euclidean normal of $\bR^n$ 
and $\mu > 1/2$.

$(b)$ the differential operators $\Lambda$ and 
the polynomials $g(x)$: 
\begin{align*}
\begin{cases}
\Lambda_i =&\frac{\p}{\p x_i}-\frac{(2\mu-1)x_i}{1-||x||^2}, \\ 
g_i(x)=& 1-||x||^2.
\end{cases}
\end{align*}
for any $1\le i\le n$.

$(c)$ the classical orthogonal polynomials  
$\{U_\alpha\,|\,\alpha \in \bN^n\}$ over the unite ball $\mathbb B^n$ in terms of $\Lambda$ 
and $g(x)$: ,
\begin{align*}
U_{\alpha}(x) =
\frac{ (-1)^{|\alpha|} (2\mu)_{|\alpha|}}{ 2^{|{\alpha}|} {|\alpha|}! 
(\mu+1/2)_{|\alpha|} }
\, \Lambda^\alpha(g^\alpha(x)),
\end{align*}
where, for any $c\in \bR$ and $k\in \bN$, 
$(c)_k=c(c+1)\cdots (c+k-1)$.

\item {\bf Classical Orthogonal Polynomials over Simplices:} 

$(a)$ $B=
T^n=\{ x\in \bR^n\, | \, \sum_{i=1}^n x_i<1; \,\, x_1, ..., x_n > 0 \}$
and the weight function
\begin{align}
w_\kappa (x)=x_1^{\kappa_1}\cdots x_n^{\kappa_n}
(1-|x|_1)^{ \kappa_{n+1}},
\end{align}
where $\kappa_i>-1$ $(1\leq i\leq n+1)$ 
and $|x|_1=\sum_{i=1}^n x_i$.

$(b)$ the differential operators $\Lambda$ and 
the polynomials $g(x)$: 
\begin{align*}
\begin{cases}
\Lambda_i &= \frac\p{\p x_i}+\frac{\kappa_i}{x_i}-
\frac{\kappa_{n+1}}{1-|x|_1},\\  
g_i(x)&= x_i(1-|x|_1)
\end{cases}
\end{align*}
for any $1\le i\le n$.

$(c)$ the classical orthogonal polynomials  
$\{U_\alpha\,|\,\alpha \in \bN^n\}$ over the simplex $T^n$ in terms of $\Lambda$ 
and $g(x)$: 
\begin{align}\label{n-Appell}
U_{\alpha}(x) =
 \Lambda^\alpha(g^\alpha(x)).
\end{align}
\end{enumerate}
\end{exam}

\begin{rmk}\label{R4.1.3}
$(a)$ A very important special family of Jacobi polynomials are the {\it Gegenbauer polynomials} which are obtained by setting $\alpha=\beta=\lambda-1/2$ for some $\lambda>-1/2$. The 
Gegenbauer polynomials are also called the 
{\it ultraspherical polynomials} in the literature.

$(b)$ For the special cases with $\lambda=0, 1, 1/2$, 
the Gegenbauer Polynomials are called the 
{\it Chebyshev polynomial of the first kind, the second kind} 
and the {\it Legendre polynomials}, respectively. 

$(c)$ When $n=2$, up to some non-zero constants the orthogonal polynomials $U_\alpha(x)$ $(\alpha\in \bN^2)$ in 
Eq.\,$(\ref{n-Appell})$ are also called {\it Appell} 
polynomials.
\end{rmk}

Note that, one important way to construct multi-variable orthogonal polynomials is to take cartesian products of orthogonal polynomials in one variable. 

More precisely, let $x=(x_1, x_2, ..., x_n)$ and 
$\{u_{i, m}(x_i)\,|\, m\ge 0\}$ $(1\le i\le n)$ be orthogonal polynomials over a subset $B_i\subset \bR$ with weight function $w_i(x_i)$ over $B_i$. Let 
\begin{align}
B&=B_1\times B_2\times \cdots \times  B_n, \\
w(x)&=w_1(x_1)w_2(x_2)\cdots w_n(x_n),\\
u_\alpha(x)&=u_{1, k_1}(x_1)u_{2, k_2}(x_2)\cdots u_{n, k_n}(x_n)
\end{align}
for any $\alpha=(k_1, k_2, ..., k_n)\in \bN^n$.

Then it is easy to see that 
$\{u_\alpha(x)\,|\, \alpha\in \bN\}$
are orthogonal polynomials over $B\subset \bR^n$ 
with respect to the weight function $w(x)$.

Furthermore, if, for any $1\le i\le n$ and $m\ge 0$, the orthogonal polynomial $u_{i, m}(x_i)=c_{i, m}\Lambda_i^m (g_i^m(x_i))$ for some nonzero $c_{i, m}\in \bR$, $g_i(x)\in \bC[x_i]$ and a differential operator $\Lambda_i$ of 
a localization of $\bC[x_i]$. 
Set $\Lambda=(\Lambda_1, \Lambda_2, ..., \Lambda_n)$ 
and $g(x)=(g_1(x_1),  g_2(x_2), ..., g_n(x_n))$.
Then, it is easy to see that $\Lambda$ is a commuting 
subset of differential operators of a localization 
of $\bC[x]$, and the orthogonal polynomials
$u_\alpha(x)$ $(\alpha\in \bN^n)$ over $B$ 
are given by
\begin{align}
u_\alpha(x)=c_\alpha \Lambda^\alpha(g^\alpha(x)),
\end{align}
where $c_\alpha=\prod_{i=1}^n c_{i, \alpha_i}$ for any 
$\alpha=(\alpha_1, \alpha_2, ..., \alpha_n)\in \bN^n$. \\
 
For later purposes, let us consider the multi-variable Jacobi orthogonal polynomials which, by Remark \ref{R4.1.3}, 
also cover the multi-variable Gegenbauer, Chebyshev and 
Legendre orthogonal polynomials. 

\begin{exam}\label{MJP}
Let $B=(-1, 1)^{\times n}\subset \bR^n$ and $\alpha, \beta\in 
\bR^n$ with all the components 
$\alpha_i, \beta_i>-1$ $(1\le i\le n)$. 
For any $1\le i\le n$, set 
\begin{align} 
\Lambda_i\!:&=\frac{d}{dx_i}-\alpha_i(1-x_i)^{-1}+\beta_i(1+x_i)^{-1},\label{MJP-Li} \\
g_i(x)&\!:=1-x_i^2.\label{MJP-gi}
\end{align}
Furthermore, set 
\begin{align}
w(x)\!:&=\prod_{i=1}^n (1-x_i)^{\alpha_i}(1+x_i)^{\beta_i},
\label{MJP-W} \\
\Lambda_{\alpha, \beta}\!:&=(\Lambda_1, 
\Lambda_2, ..., \Lambda_n),\label{MJP-Lambda}\\
g(x)\!:&=(g_1(x), g_2(x), ..., g_n(x)). \label{MJP-g}
\end{align}

For any $1\le i\le n$ and $m\ge 0$, let 
$P^{\alpha_i, \beta_i}_m (x)$ be the $m^{th}$ one-variable Jacobi polynomial in $x_i$ $($See Example \ref{COP}, $(3)$$)$ with 
$\alpha=\alpha_i$ and $\beta=\beta_i$. 
For any ${\bf m}=(m_1, m_2, ..., m_n)\in \bN^n$, set 
\begin{align}
P^{\alpha, \beta}_{\bf m} (x)\!:=\prod_{i=1}^n 
P^{\alpha_i, \beta_i}_{m_i}(x_i). \label{MJP-P}
\end{align}

Then, for any fixed $\alpha, \beta\in (\bR^{>-1})^{\times n}$,
the sequence $\{P^{\alpha, \beta}_{\bf m} (x)\,|\, {\bf m}
\in \bN^n\}$ forms a sequence  of orthogonal polynomials over 
$B$ with the weight function given by Eq.\,$(\ref{MJP-W})$.
From Eq.\,$(\ref{JP-LP})$, it is easy to see that the relation 
of $\{P^{\alpha, \beta}_{\bf m} (x)\,|\, {\bf m}
\in \bN^n\}$ with the commuting differential operators 
$\Lambda$ in Eq.\,$(\ref{MJP-Lambda})$ and the polynomial $g(x)$ in Eq.\,$(\ref{MJP-g})$ is given by 
\begin{align} \label{MJP-LP}
P^{\alpha, \beta}_{\bf m} (x)=\frac {(-1)^{|{\bf m}|}}
{2^{|\bf m|}{\bf m!} }\, \Lambda_{\alpha, \beta}^{\bf m} 
g^{\bf m}(x).
\end{align}
\end{exam}

Finally, let us summarize the relations of orthogonal polynomials 
with commuting differential operators of order one with constant leading coefficients in the following lemma.

\begin{lemma}\label{L4.1.5}
$(a)$ Up to some nonzero multiplicative scalars, all classical orthogonal polynomials $\{u_\alpha(x)\,|\, \alpha\in \bN^n\}$ 
above including those obtained by Cartesian products of 
classical orthogonal polynomials have the form 
in Eq.\,$(\ref{LaGa})$ for some  
$g(x)=(g_1(x_1),  g_2(x_2), ..., g_n(x_n))\in \bC[x]^n$ and     differential operators $\Lambda=(\Lambda_1, \Lambda_2, ..., \Lambda_n)$ of some localizations $\cB$ of $\bC[x]$. 

$(b)$ The set $\Lambda$ is a commuting subset of differential operators of $\cB$ of order one with constant leading coefficients.

$(c)$ For any nonzero $\alpha\in \bN^n$, the orthogonal polynomials $u_\alpha(x)\in \im' \Lambda\!:=\bC[x]\cap\sum_{i=1}^n (\Lambda_i \bC[x])$.
\end{lemma}

Note that $(a)$ and $(b)$ follow immediately from the discussion in this subsection. $(c)$ follows from the fact that 
$\Lambda^\alpha (g_i(x)g^\alpha)$ for any $1\le i\le n$ and 
$\alpha\in \bN^n$, which can also be easily checked directly. 

\subsection{The Image Conjecture for the Differential Operators Related with the Multi-Variable Jacobi Orthogonal Polynomials} \label{S4.2} 

Considering the important roles of orthogonal polynomials 
played in so many different areas, it will be interesting 
to see if the {\it image} conjecture ({\bf IC}), Conjecture \ref{IC}, 
also holds for the commuting differential operators related with 
orthogonal polynomials.

In this subsection, we consider the straightforward generalization of the \IC  for the following family of commuting differential operators of Laurent polynomial algebras in a slightly more general setting, namely, with the base field $\bC$ replaced by integral domains over $\bC$. As we will see that the straightforward generalization of the \IC is false for these differential operators (See Propositions \ref{P4.2.1} and \ref{P4.2.2}). But another generalization of the \IC for these  differential operators under the constraints from the multi-variable Jacobi orthogonal polynomials actually holds (See Corollary \ref{C4.2.3}).

Let $\cA$ be any  integral domain over 
$\bC$ and $\cAzz$ the algebra of Laurent polynomials 
with coefficients in $\cA$. 
For any $\lambda=(\lambda_1, \lambda_2, ..., \lambda_n)\in 
\bC^n$,  we set $\Phi_{\lambda_i}\!:=\p_i+\lambda_i z_i^{-1}$ 
$(1\le i\le n)$ and $\Phi_\lambda\!:=(\Phi_{\lambda_1}, \Phi_{\lambda_2}, ..., \Phi_{\lambda_n})$. We will also view $\Phi_{\lambda}$ as a commuting subset (instead of just an $n$-tuple) 
of differential operators of $\cAzz$ 
of order one with constant leading coefficients.

Note that, the differential operators 
$\Phi_\lambda$ are essentially the 
differential operators related with 
the multi-variable Jacobi orthogonal 
polynomials $P_{\bf m}^{\alpha, \beta}(x)$ 
$({\bf m}\in \bN^n)$ in Eq.\,(\ref{MJP-P})  with 
$\alpha=0$ or $\beta=0$.
For example, by setting $\alpha=0$ and 
$\beta=\lambda$, and changing the variables 
$x_i\to z_i-1$ $(1\le i\le n)$, 
from Eqs.\,(\ref{MJP-Li})  and (\ref{MJP-Lambda}) 
we see that the differential 
operators $\Lambda_{\alpha, \beta}$ 
related with the Jacobi polynomials will 
coincide with the differential operators 
$\Phi_\lambda$. Similarly, this is also 
the case when $\beta=0$ if we set 
$\alpha=\lambda$ and apply the changing 
of variables $x_i\to z_i+1$ $(1\le i\le n)$.

But we emphasize that, unlike for the parameters 
$\alpha$ and $\beta$ of the Jacobi polynomials, 
here we do not require $\lambda_i>-1$ $(1\le i\le n)$ 
nor even $\lambda\in \bR^n$ unless stated otherwise.

Now we fix any $\lambda\in\bC^n$ and 
the differential operators $\Phi_\lambda$ as above, 
and set $\im\Phi_\lambda\!:=\sum_{i=1}^n 
(\Phi_{\lambda_i} \cAzz)$. We also fix the following notation 
that will be used throughout the rest of this paper.

For any $\gamma\in \bZ^n$ and $g(z)\in\cAzz$, we denote by 
$[z^\gamma]g(z)$ the coefficient of the monomial $z^\gamma$ 
in $g(z)$. For convenience,  
we also allow $\gamma$ in the notation above 
to be any element of $\bC^n$, i.e. 
we set $[z^\gamma]g(z)=0$ for any $g(z)\in\cAzz$ and 
$\gamma\in \bC^n\backslash \bZ^n$. In the case that 
$[z^\gamma]g(z)=0$, we also say that $g(z)$ has 
no $z^\gamma$ term.

With all the notations fixed above, 
we have the following proposition.

\begin{propo}\label{P4.2.1}
For any $\lambda\in \bC^n$, denote by the abusing notation 
$-\lambda-1$ the $n$-tuple $-\lambda-(1, 1, ..., 1)$. Then, 
we have

$(a)$ $\im\Phi_\lambda$ is the $\cA$-subspace of $\cAzz$ consisting 
of the Laurent polynomials $g(z)\in \cAzz$ 
with $[z^{-\lambda-1}]g(z)=0$.  In particular,
$\im\Phi_\lambda=\cAzz$ if $\lambda \not\in \bZ^n$. 

$(b)$ $\im \Phi_\lambda$ is a Mathieu subspace of $\cAzz$ 
iff $\lambda \not\in \bZ^n$ or $\lambda=(-1, -1, ...,-1)$.
\end{propo} 

\pf We first prove the proposition for the case $n=1$, 
i.e. for the one-variable case.

For any $g(z)\in \cAzz$, consider 
the ordinary differential 
equation with the unknown function 
$f(z)\in \cAzz$:
\begin{align}\label{P4.2.1-pe1}
\Phi_\lambda f = f'+\lambda z^{-1}f=g.
\end{align}

The equation above can be solved by the following standard trick  in ODE. First, we view the equation as 
a differential equation for elements of 
$\cA[z^{\pm\lambda}, z^{\pm 1}]$,  
and set $\tilde f(z)\!:=z^{\lambda}f(z)\in 
\cA[z^{\pm\lambda}, z^{\pm 1}]$. 
Then $f(z)=z^{-\lambda}\tilde f(z)$. Plug this expression of 
$f(z)$ in Eq.\,(\ref{P4.2.1-pe1}), it is easy to check that 
$\tilde f(z)$ satisfies the following equation:
\begin{align*}
z^{-\lambda}\tilde f' = g.
\end{align*}

Therefore, we have 
$\tilde f(z) = \int z^{\lambda} g(z)\, dz$  
and
\begin{align}
f(z)&= z^{-\lambda}\int z^{\lambda} g(z)\, dz. \label{P4.2.1-pe5}
\end{align}

From the arguments, we see that any solution 
$f\in \cAzz$ of Eq.\,(\ref{P4.2.1-pe1}) must 
be given by Eq.\,(\ref{P4.2.1-pe5}) up to 
an $z^{-\lambda}$ term. 
But, conversely, the RHS of 
Eq.\,(\ref{P4.2.1-pe5}) 
does not necessarily produce an element of 
$\cAzz$ unless $z^\lambda g(z)$ has 
no $z^{-1}$ term, i.e. the residue 
$\res\, z^\lambda g(z)=0$.

Therefore, the differential equation Eq.\,(\ref{P4.2.1-pe1}) has a Laurent polynomial solution $f(z)\in \cAzz$ iff $\res\, z^\lambda g(z)=0$ iff  
$[z^{-\lambda-1}]g(z)=0$. Hence, we have $(a)$ 
of the proposition for the case $n=1$.

To show $(b)$ for the case $n=1$, let us look at all the values of 
$\lambda\in \bC$ such that $\im \Phi_\lambda$ is 
a Mathieu subspace of $\cAzz$. 

First, if $\lambda\not \in \bZ$, by $(a)$ 
$\im \Phi_\lambda=\cAzz$ which is obviously 
a Mathieu subspace of $\cAzz$. 

Consider the case $\lambda\in \bZ$. 
If $-\lambda-1\neq 0$, i.e. $\lambda\neq -1$.
Then, by statement $(a)$ for the case $n=1$, 
we have, $1\in \im \Phi_\lambda$ and 
$\im \Phi_\lambda \ne \cAzz$. 
By the general property of Mathieu subspaces given 
in Lemma \ref{No-1-lemma} in Section \ref{S5},   
$\im \Phi_\lambda$ is not a Mathieu subspace 
of $\cAzz$. 

A more convincing counter-example for this case can 
be constructed as follows. Set 
\allowdisplaybreaks{
\begin{align}
v(z)\!:&=z^{-\lambda-1} \label{P4.2.1-pe3} \\
u(z)\!:&=\begin{cases} 1+z^{-\lambda} &\mbox{ if } \lambda<-1; \\ 
1+z^{-\lambda-2} &\mbox{ if } \lambda>-1. 
\end{cases}\label{P4.2.1-pe4}
\end{align} } 
Then, it is easy to check that,  
for any $m\ge 1$, we have $[z^{-\lambda-1}]u^m=0$ 
and $[z^{-\lambda-1}](u^mv)=z^{-\lambda-1}$. 
By statement $(a)$ for the case $n=1$, 
we have that, for any $m\ge 1$, 
$u^m \in \im \Phi_\lambda$ but 
$u^m v\not \in \im \Phi_\lambda$. 

Next, consider the case $\lambda=-1$. By $(a)$, 
we know that $\im \Phi_{\lambda=-1}$ is the $\cA$-subspace 
of $\cAzz$ consisting of Laurent polynomials 
with no constant terms. In the case that 
$\cA=\bC$, $(b)$ follows directly from 
Duistermaat and van der Kallen's theorem, 
Theorem \ref{ThmDK}. 

In the case that $\cA \neq \bC$, $(b)$ also follows from 
Theorem \ref{ThmDK} via Lefschetz's principle since,  
whenever we fix $a(z), b(z)\in \cAzz$ with   
$a^m(z)$ has no constant term for any $m\ge 1$,  
to show that $a^m b$ 
has no constant term when $m>>0$, 
we may replace $\cA$ by the field $K$ generated by 
the (finitely many) coefficients of $a(z)$ and $b(z)$ 
over $\bQ$, which can be embedded in $\bC$ 
as a subfield.

Therefore, we have proved the proposition for the 
case $n=1$. Now we assume $n\ge 2$.

For convenience, throughout the rest of the proof, 
we denote by ${\mathcal U}$ the subalgebra 
of $\cAzz$ of Laurent polynomials in $z_i$ 
$(2\le i\le n)$ with coefficients in $\cA$. 
Note that $\cAzz$ may be viewed as 
the Laurent polynomial algebra in $z_1$ 
over $\cU$, i.e. we have $\cAzz=\cU[z_1^{-1}, z_1]$.
 
First, assume that $\lambda\not \in \bZ^n$. 
Then there exists $1\le i\le n$ such that 
$\lambda_i \not \in \bZ$. Without losing any generality, 
we assume $\lambda_1\not \in \bZ$.  
By statement $(a)$ for the one-variable case with 
$\cA$ replaced by $\cU$ and $\lambda$ by $\lambda_1$, 
we have 
\begin{align*}
\Phi_{\lambda_1}(\cU[z_1^{-1}, z_1])=
\cU[z_1^{-1}, z_1]=\cAzz.
\end{align*}

Since $\Phi_{\lambda_1}(\cU[z_1^{-1}, z_1])=
\Phi_{\lambda_1}(\cAzz) \subset \im \Phi_{\lambda}$, 
from the equation above, we have 
$\im \Phi_{\lambda}=\cAzz$. 

Now assume $\lambda \in \bZ^n$. We first show that   
$z^\beta \in \im\Phi_\lambda$ for any 
$\beta\in \bZ^n$ with $\beta\ne -\lambda-1$.

Pick up any $\beta \ne -\lambda-1\in \bZ^n$, 
there exists $1\le i\le n$ such that the $i^{th}$ 
component $\beta_i$ of $\beta$ is different 
from $-\lambda_i-1$. We assume 
$\beta_1 \neq -\lambda_1-1$ (the proof for other cases 
is similar). Then, by statement $(a)$ for the one-variable case with 
$\cA$ replaced by $\cU$ and $\lambda$ by $\lambda_1$, 
we have 
\begin{align*}
z^\beta =(z_2^{\beta_2} z_3^{\beta_3}\cdots z_n^{\beta_n}) z_1^{\beta_1}  \in \Phi_{\lambda_1}(\cU[z_1^{-1}, z_1])=
\Phi_{\lambda_1}(\cAzz) \subset \im\Phi_\lambda.
\end{align*}

Consequently, we see 
that any $g(z)\in \cAzz$ with 
$[z^{-\lambda-1}]g(z)=0$ lies in 
$\im \Phi_\lambda$. 
Conversely, for any $g(z)\in \im\Phi_\lambda$, 
we write $g(z)$ as
\begin{align}\label{P4.2.1-pe7}
g(z)=\sum_{i=1}^n \Phi_{\lambda_i} f_i(z)
\end{align}  
for some $f_i(z)\in \cAzz$ $(1\le i\le n)$.

By statement $(a)$ for the one-variable case with 
$\cA$ replaced by $\cU$ and $\lambda$ by $\lambda_1$, 
we see that $\Phi_{\lambda_1} f_1(z)$ 
in Eq.\,(\ref{P4.2.1-pe7}) has no 
$z^{\gamma}$ term for any $\gamma\in\bZ^n$ 
with the first component $\gamma_1=-\lambda_1-1$. 
Hence, for the similar reason, 
for any $2\le i \le n$, 
$\Phi_{\lambda_i} f_i(z)$ in Eq.\,(\ref{P4.2.1-pe7}) 
can not have the $z^\gamma$ term for 
any $\gamma\in\bZ^n$ with the 
$i^{th}$ component $\gamma_i=-\lambda_i-1$. 
Therefore, by Eq.\,(\ref{P4.2.1-pe7}) 
we see that $g(z)$ can not have any  
$z^{-\lambda-1}$ term. 

Combining the results in the last two paragraphs,  
we have statement $(a)$ of the proposition. 

Next we show statement $(b)$. First, 
if $\lambda \not \in \bZ^n$, by $(a)$ we have  
$\im \Phi_\lambda=\cAzz$ which is obviously a Mathieu 
subspace of $\cAzz$.

Assume that $\lambda\in \bZ^n$ but 
$\lambda\neq (-1, -1, ..., -1)$.
Since $-\lambda-1\ne 0\in \bN^n$, by statement $(a)$, 
we have, $1\in \im \Phi_\lambda$ and 
$\im \Phi_\lambda \ne \cAzz$. 
By Lemma \ref{No-1-lemma} in Section \ref{S5}, 
$\im \Phi_\lambda$ is not a Mathieu subspace 
of $\cAzz$.

Finally, consider the case $\lambda=(-1, -1, ..., -1)$. 
Similarly as for the one variable case, by Duistermaat and van der Kallen's theorem, Theorem \ref{ThmDK}, and Lefschetz's principle, it is easy to see that $\im \Phi_\lambda$ is a Mathieu subspace of $\cAzz$ in this case. 
\epfv 

\begin{rmk}
From Proposition \ref{P4.2.1}, we see that 
$($the straightforward generalization of$)$ 
the \IC for the Laurent polynomial algebras 
does not always hold.

But, on the other hand, it is still interesting to see, 
for which commuting differential operators of $\czz$ or any localization of $\cz$, the \IC does hold. For example,   
for the differential operators $\Phi_\lambda$ with 
$\lambda=(-1, -1, ..., -1)$, the \IC is 
actually equivalent to Duistermaat and 
van der Kallen's theorem, Theorem \ref{ThmDK}.
Even more mysteriously, among all the cases 
that $\im \Phi_\lambda \neq \cAzz$, 
this is the only case that the \IC holds.
\end{rmk}

Next let us consider the \IC for the polynomial algebra 
$\cA[z]$ (instead of $\cAzz$) 
and the differential operators $\Phi_\lambda$ 
(even though $\cA[z]$ is not closed 
under the action of $\Phi_\lambda$). 
But, first, we need prove the 
following lemma.

\begin{lemma}\label{P4.2.2-L}
Let $\cA$, $\Phi_\lambda$ $(\lambda\in \bC^n)$ 
and $\im\Phi_\lambda$ as 
in Proposition \ref{P4.2.1}. Then, for any $g(z)\in \cAz$,  
$g(z)\in \im \Phi_\lambda$ iff 
$g(z)\in \sum_{i=1}^n 
(\Phi_{\lambda_i}\cA[z])$.
\end{lemma}

\pf Note that the $(\Leftarrow)$ part
of the lemma is trivial. We use induction on 
$n\ge 1$ to show the other part.

First, assume $n=1$ and $g(z)\in \im \Phi_\lambda$. 
Then, by Proposition \ref{P4.2.1}, $(a)$, 
we have $[z^{-\lambda-1}]g(z)=0$. 
By writing $g(z)$ as a linear combination 
of monomials $z^k$ $(k\in \bN)$ over $\cA$, 
it is easy to check that  
Eq.\,$(\ref{P4.2.1-pe5})$ 
has a polynomial solution $f(z)\in \cAz$. 
Since $\Phi_\lambda f(z)=g(z)$ as shown 
in the proof of Proposition \ref{P4.2.1}, 
the lemma holds in this case.

Next assume the lemma holds for the $n-1$ case and 
consider the $n$-variable case.

Note first that, for any $\lambda\in \bC^n$ 
and $g(z)\in \cAz$, 
there exist $u(z) \in \cA[z_2, z_3, ..., z_n]$
and $v(z)\in \cAz$ such that $[z^\gamma]v(z)=0$ 
for any $\gamma\in \bN^n$ with the first component 
$\gamma_1=-\lambda_1-1$ and 
\begin{align}\label{P4.2.2-L-pe1}
g(z)=z_1^{-\lambda_1-1}u(z)+v(z).
\end{align}

Now further assume $g(z)\in \im \Phi_\lambda$. Then,
by Proposition \ref{P4.2.1}, $(a)$, we have $[z^{-\lambda-1}]g(z)=0$.
We show below that both terms on the right hand side of 
Eq.\,(\ref{P4.2.2-L-pe1}) lie in 
$\sum_{i=1}^n (\Phi_{\lambda_i} \cA[z])$.

First, if $-\lambda_1-1\not \in \bN$, we have $u(z)=0$. 
Otherwise, set $\cA'\!:=\cA[z_1]$ and 
$z''\!:=(z_2, ..., z_n)$. 
We view $u(z)$ and also $z_1^{-\lambda_1-1}u(z)$ 
as polynomials in $z''$ over the integral domain 
$\cA'$. Note that, the coefficient 
of the monomial $z_2^{-\lambda_2-1}\cdots z_n^{-\lambda_n-1}$ in 
$u(z)\in \cA'[z'']$ is same as $[z^{-\lambda-1}]g(z)$ 
which is equal to zero. Hence the coefficient 
of $z_2^{-\lambda_2-1}\cdots z_n^{-\lambda_n-1}$ in  
$z_1^{-\lambda_1-1} u(z)\in  \cA'[z'']$ is also 
equal to zero. 

Apply Proposition \ref{P4.2.1}, $(a)$ to 
$z_1^{-\lambda_1-1} u(z) \in  \cA'[z'']$ 
with $\cA$ replaced by $\cA'$, we know that 
$z_1^{-\lambda_1-1} u(z)\in \sum_{i=2}^n (\Phi_{\lambda_i} 
\cA'[z_2^\pm, ..., z_n^\pm])$. Then, 
by applying the induction assumption 
to $z_1^{-\lambda_1-1}u(z)$ with $\cA$ 
replaced by $\cA'$, we have,  
$z_1^{-\lambda_1-1}u(z)\in 
\sum_{i=2}^n (\Phi_{\lambda_i} \cA'[z_2, ... , z_n])
=\sum_{i=2}^n (\Phi_{\lambda_i}\cAz)$.

Second, set $\cA''\!:=\cA[z_2, ... , z_n]$. Then, 
viewing $v(z)$ as a polynomial in $z_1$ 
over the integral domain $\cA''$, we have 
$[z_1^{-\lambda_1-1}]v(z)=0$. 
By Proposition \ref{P4.2.1}, $(a)$ with $\cA$ replaced by  
$\cA''$,  
$v(z)\in \Phi_{\lambda_1} (\cA''[z_1^{-1}, z_1])$.
Applying the lemma for the case $n=1$ to $v(z)$ 
with $\cA$ replaced by $\cA''$,  
we  have $v(z)\in \Phi_{\lambda_1}(\cA''[z_1]) 
=\Phi_{\lambda_1}(\cA[z])$. 
Then, by Eq.\,(\ref{P4.2.2-L-pe1}), 
we have $g(z)\in \sum_{i=1}^n 
(\Phi_{\lambda_i} \cA[z])$,  
and the lemma follows.
\epfv

\begin{propo}\label{P4.2.2}
Let $\cA$ and $\Phi_\lambda$ $(\lambda\in \bC^n)$ as 
in Proposition \ref{P4.2.1}. Set 
\begin{align}\label{P4.2.2-e1}
\im'\Phi_\lambda\!:=\cAz \cap \sum_{i=1}^n
(\Phi_{\lambda_i} \cA[z]). 
\end{align}
Then,  
$(a)$ $\im'\Phi_\lambda=\cA[z]$ iff $\lambda \not \in 
(\bZ^{<0})^n$, where $\bZ^{<0}$ denotes the set 
of all negative integers.

$(b)$ $\im'\Phi_\lambda$ is a Mathieu subspace of $\cA[z]$ iff 
$\lambda \not \in (\bZ^{<0})^n$ or $\lambda=(-1, -1, ..., -1)$. 
\end{propo}

\pf  Note first that, by Eq.\,(\ref{P4.2.2-e1}) and Lemma \ref{P4.2.2-L}, it is easy to see that 
\begin{align}\label{C4.2.2-pe1}
\im'\Phi_\lambda=\cA[z]\cap\im\Phi_\lambda.
\end{align}

If $\lambda\not \in \bZ^n$, by Proposition \ref{P4.2.1}, 
$(a)$ and the equation above, we have that 
$\im'\Phi_\lambda=\cA[z]$.

Assume $\lambda\in \bZ^n$ and set 
$\alpha\!:=-\lambda-(1, 1, ..., 1)$.
By Proposition \ref{P4.2.1}, 
$(a)$ and Eq.\,(\ref{C4.2.2-pe1}), 
we know that $\im'\Phi_\lambda$ is the $\cA$-subspace 
of polynomials with no $z^\alpha$ term. 

If $\lambda \not \in (\bZ^{<0})^n$, 
i.e. $\lambda_i\ge 0$ for some $1\le i\le n$, then 
we have $\alpha\not \in \bN^n$ and hence 
$\im'\Phi_\lambda=\cA[z]$. 

Now assume $\lambda \in (\bZ^{<0})^n$ but 
$\lambda\neq (-1, -1, ..., -1)$. Then there exists 
$1\le j\le n$ such that $\lambda_j \le -2$.
Note that, in this case $\alpha\in \bN^n$ 
but $\alpha\ne 0$. Consequently, we have,  
$1\in \im'\Phi_{\lambda}$ and 
$\im'\Phi_{\lambda}\ne \cAz$. Then, 
by Lemma \ref{No-1-lemma} 
in Section \ref{S5}, $\im'\Phi_{\lambda}$ 
is not a Mathieu subspace of $\cAz$.

Finally, if $\lambda=(-1, -1, ..., -1)$, then $\alpha=0$ 
and $\im'\Phi_\lambda$ is the ideal of all polynomials 
with no constant terms. So in this case 
$\im'\Phi_{\lambda}\neq \cAz$ but is the ideal of $\cAz$ 
generated by $z_i$ $(1\le i\le n)$. Hence it is 
a Mathieu subspace of $\cAz$.

So we have exhausted all possible choices of 
$\lambda\in \bC^n$. Combining all the results above, 
it is easy to see that both $(a)$ and $(b)$ 
of the proposition hold. 
\epfv 

As pointed out at the beginning of this subsection, 
up to some changes of variables, the 
differential operators $\Phi_\lambda$ are 
same as the differential operators 
$\Lambda_{\alpha, \beta}$ in 
Eq.\,(\ref{MJP-Lambda}) related with 
the multi-variable Jacobi orthogonal 
polynomials $P^{\alpha, \beta}_{\bf m}$ 
$({\bf m}\in \bN^n)$ in Eq.\,(\ref{MJP-P}) 
with $\alpha=\lambda$ and $\beta=0$ 
or $\alpha=0$ and $\beta=\lambda$. 

Note that, the constraints on 
the parameters $\alpha$ and $\beta$ of 
the Jacobi polynomials are 
$\alpha, \beta\in \bR^n$ and 
the components $\alpha_i, \beta_i>-1$ 
$(1\le i\le n)$. 
Now, if we put the same 
constraints on $\lambda$, then, by 
Proposition \ref{P4.2.2},
it is easy to see that we have 
the following corollary.

\begin{corol}\label{C4.2.3}
Let $\cA$ and $\Phi_\lambda$ as 
in Proposition \ref{P4.2.1}. 
Assume further that $\lambda\in \bR^n$ 
and $\lambda_i>-1$ for any $1\le i\le n$. 
Then, $\im'\Phi_\lambda=\cAz$ and hence is a 
Mathieu subspace of $\cAz$.
\end{corol}

\renewcommand{\theequation}{\thesection.\arabic{equation}}
\renewcommand{\therema}{\thesection.\arabic{rema}}
\setcounter{equation}{0}
\setcounter{rema}{0}

\section{\bf Generalizations of the Image Conjecture}\label{S4b}

Motivated by the discussions in Subsection \ref{S4.2}, 
we first formulate in Subsection \ref{S4.3}  
a generalization (See Conjecture \ref{GIC}) 
of the {\bf IC} for the 
commuting differential operators 
related with orthogonal polynomials. 
We also show that, beside a trivial case, 
Conjecture \ref{GIC} is actually 
equivalent to a special case of another conjecture, 
Conjecture \ref{GMC}, on integrals of polynomials 
over open subsets of $\bR^n$ with any  
positive measures. As we will see that  
Conjecture \ref{GMC} can also be viewed as 
a natural variation of the Mathieu conjecture, 
Conjecture \ref{MC}, and Duistermaat and 
van der Kallen's theorem, Theorem \ref{ThmDK}.

In Subsection \ref{S4.4}, we will prove some cases of 
Conjectures \ref{GIC} and \ref{GMC}. We will also 
discuss a connection of Conjecture \ref{GMC} 
with the so-called {\it polynomial moment problem}.

\subsection{The Generalized Image Conjecture}\label{S4.3}

First, we propose the following generalization of the {\bf IC} for the commuting differential operators related with classical 
orthogonal polynomials. 

\begin{conj}\label{GIC}
Let $B\subset \bR^n$, $w(z)$ and 
$\{u_\alpha\,|\, \alpha\in \bN^n\}$ as in Definition 
\ref{Def-OP} with $x$ replaced by $z$. Assume further 
that the orthogonal polynomials 
$\{u_\alpha\,|\, \alpha\in \bN^n\}$ 
can be obtained via Eq.\,$(\ref{LaGa})$ 
for some commuting differential operators 
$\Lambda=(\Lambda_1, \Lambda_2, ..., \Lambda_n)$ of a localization of the polynomial algebra $\cz$. 
Set 
\begin{align}\label{Def-im'}
\im'\Lambda\!:=\cz\,\bigcap \, \sum_{i=1}^n (\Lambda_i \cz). 
\end{align}
Then, $\im'\Lambda$ is a Mathieu subspace of $\cz$.
\end{conj}

Note first that, when $\Lambda$ are differential 
operators of $\cz$ (instead of a localization of $\cz$), 
we have, $\im'\Lambda=\im\Lambda$. Therefore, the conjecture above can be viewed as a generalization of the {\it image} conjecture, Conjecture \ref{IC}, to the differential operators related with classical orthogonal polynomials.

Note also that, by Lemma \ref{L4.1.5}, $(c)$,  
$\im'\Lambda$ has co-dimension in $\cz$ zero or one depending 
on whether $u_0(z)$ lies in $\im'\Lambda$ or not. 
Since $u_0(z)$ is a nonzero constant, we have 
$\im'\Lambda=\cz$ iff $1\in \mbox{Im}'\Lambda$.
 
Therefore, if $1\in \mbox{Im}'\Lambda$, we have $\im'\Lambda=\cz$, 
and hence Conjecture \ref{GIC} holds trivially in this case. 
If $1\not \in \im'\Lambda$, we have 
\begin{align}\label{Im=Su}
\mbox{Im}'\Lambda=\mbox{Span}_\bC\{u_\alpha(z)\,|\, \alpha\ne 0\}.
\end{align}

In this case, Conjecture \ref{GIC} turns 
out to be equivalent to a special case of 
the following conjecture.

\begin{conj}\label{GMC} 
Let $B$ be any non-empty open subset of $\bR^n$ and   
$\sigma$ any positive 
measure such that $\int_B g(z)\, d\sigma$ is finite for 
any $g(z)\in \cz$. Let $\cM_B(\sigma) $ be the subspace 
of all polynomials $f(x)\in \cz$ such that 
$\int_B f(x) \, d \sigma =0$.
Then $\cM_B(\sigma)$ is a Mathieu subspace of $\cz$.
\end{conj}

\begin{propo}\label{P4.2.3}
Let $B\subset \bR^n$, $w(z)$ and $\Lambda$ be 
as in Conjecture \ref{GIC}. Let $\sigma$ be the measure 
on $B$ such that $d\sigma=w(z)dz$. 
Assume further that $1\not \in \im'\Lambda$. 
Then, we have

$(a)$ $\im'\Lambda=\cM_B(\sigma)$.  

$(b)$ Conjecture \ref{GIC} for the differential operators 
$\Lambda$ is equivalent 
to Conjecture \ref{GMC} for the open subset $B\subset \bR^n$ with the  measure $\sigma$. 
\end{propo}

\pf First it is easy to see that $(b)$ follows directly from $(a)$. 

To show $(a)$, choose any $f(z)\in \cz$ and write it (uniquely) 
as $f(z)=\sum_{\alpha\in \bN^n}c_\alpha u_\alpha(z)$ with 
$c_\alpha\in \bC$. Then, by Eq.\,(\ref{Im=Su}) and the assumption that 
$1\not \in \im'\Lambda$, we have that, $f(z)\in \im'\Lambda$
iff $c_0=0$. 

On the other hand, since $\{u_\alpha\,|\, \alpha\in \bN^n\}$ is  
an orthogonal basis of $\cz$ with respect to 
the Hermitian form defined in Definition \ref{Def-OP}, 
for any $\alpha\neq 0$, we have
\begin{align} 
\int_B u_\alpha(z) w(z)dz=
\bar u_0^{-1} \int_B  u_\alpha(z) \bar u_0 \,w(z)dz=0.
\end{align}
Therefore, we have
\begin{align*}
\int_B f(z) \, d\sigma=
\int_B f(z) w(z)dz=c_0\int u_0(z) w(z)dz =
c_0 u_0\int_B w(z)\, dz.  
\end{align*}

Since $u_0$ and $\int w(z)\, dz$ are nonzero constants, 
we have that, $f(z)\in \cM_B(\sigma)$ 
iff $c_0=0$.

Combining the results above, we have 
$\im'\Lambda=\cM_B(\sigma)$ which is 
$(a)$ of the proposition.
\epfv

\begin{rmk}
For the special case that $1\in \im'\Lambda$, 
as pointed out above Conjecture \ref{GIC} holds trivially. 
But in this case, Conjecture \ref{GMC}, even for the 
measures $d\sigma=w(z)dz$ given by the weight functions $w(z)$  
of classical orthogonal polynomials, can still be highly
non-trivial. See Subsection \ref{S4.4} 
for more discussions. 
\end{rmk}

\vskip1mm

Several more remarks on Conjecture \ref{GMC} are as follows.

\vskip2mm

First, Conjecture \ref{GMC} can be viewed as a 
variation of the Mathieu conjecture, Conjecture 
\ref{MC}, with the reductive Lie group $G$ replaced by 
an open subset $B\subset \bR^n$; the Haar measure 
by any positive measure $\sigma$; and $G$-finite 
functions by polynomials. 
Furthermore, by Eq.\,(\ref{Im=Su}) and  
Proposition \ref{P4.2.3}, $(a)$, we see that  
Conjecture \ref{GMC} with $d\sigma=w(z)dz$ 
can also be viewed as a natural variation 
of Duistermaat and van der Kallen's theorem, 
Theorem \ref{ThmDK}, 
with the basis of $\cz$ formed by monomials 
$z^\alpha$ $(\alpha\in \bN^n)$ replaced by 
the basis formed by the orthogonal polynomials 
$u_{\alpha}(z)$ $(\alpha\in \bN^n)$. 
More precisely, Conjecture \ref{GMC} 
with $d\sigma=w(z)dz$ can be 
re-stated as follows.

\begin{conj} \label{GDK}
Let $\{u_\alpha\,|\, \alpha\in \bN^n \}$ be a 
sequence of orthogonal polynomials over $B\subset \bR^n$ 
with the weight function $w(z)$. 
Let $\cM$ be the subspace of $f(z)\in \cz$ whose  
\underline{constant term} $($the coefficient of $u_0$$)$ 
in the unique expansion of $f(z)$ in terms of 
$u_\alpha(z)$'s is equal to zero. 
Then $\cM$ is a Mathieu subspace of 
the polynomial algebra $\cz$.
\end{conj}

Of course, a shorter way to state the conjecture above 
is that the subspace $\cM$ spanned by the orthogonal 
polynomials $u_\alpha(z)$ $(\alpha\ne 0)$ over $\bC$ 
is a Mathieu subspace of the polynomial algebra $\cz$. 

The second remark is that Conjecture \ref{GMC} in general does not hold 
for analytic functions.

\begin{exam}\label{E4.2.5}
Let $B=(0, \, 1)\subset \bR$ and $d\sigma=dz$. Let 
$f(x)=e^{2\pi\sqrt{-1}\, z}$ and $g(z)=z$. 
Then, for any $m\ge 1$, we have 
\begin{align*}
\int_0^1 f^m(x) \, d\sigma=
\int_0^1 e^{2m\pi\sqrt{-1}\, z} \, dz =0. 
\end{align*}
But  
\begin{align*}
\int_0^1 f^m(z) g(z)\, d\sigma&=
\int_0^1 z e^{2m\pi\sqrt{-1}\, z} \, dz \\
&=\frac 1{2m\pi \sqrt{-1}}\left(\left. ze^{2m\pi\sqrt{-1}\, z}\right|_0^{1}- \int_0^1 e^{2m\pi\sqrt{-1}\, z} \, dz \right)\\
&=\frac1{2m\pi\sqrt{-1}}\ne 0. 
\end{align*}
\end{exam}

The third remark is that Conjecture \ref{GMC} does not hold 
without the positivity assumption on the 
measure $\sigma$.

\begin{exam}\label{E4.2.6}
Let $B=(-1, 1)\subset \bR$ and $d\sigma=zdz$. Let 
$f(z)=z^2$ and $g(z)=z$. 
Then, for any $m\ge 1$, we have 
\begin{align*}
\int_{-1}^1 f^m(z) \, d\sigma=
\int_{-1}^1 z^{2m+1} \, dz =0 
\end{align*}
But
\begin{align*} 
\int_{-1}^1 f^m(z) g(z)\, d\sigma=
\int_{-1}^1 z^{2m+2}\, dz =\frac2{2m+3}\ne 0 
\end{align*}
\end{exam}

Finally, one interesting observation about the example 
above is as follows. Even though Conjecture \ref{GMC} 
fails for this example, if we consider the 
differential operator $\Lambda$ related with 
the ``weight" function $w(z)=z$ as 
in Eq.\,(\ref{L4OP}), namely,
\begin{align*} 
\Lambda=\frac d{dz}+w^{-1}(z)\frac {dw(z)}{dz} 
=\frac d{dz}+z^{-1},
\end{align*}
then, by Corollary \ref{C4.2.3}, 
we see that 
Conjecture \ref{GIC} (formally)   
for the differential operator 
$\Lambda$ above still holds.

\subsection{Some Cases of Conjectures \ref{GIC} and \ref{GMC}}\label{S4.4} 

Despite the simple appearances of Conjectures 
\ref{GIC} and \ref{GMC}, there are only few cases that are known  
for these two conjectures. 

First, for the differential operators related with the multi-variable Jacobi orthogonal polynomials $($See Example \ref{MJP}$)$, we have the 
following proposition. 

\begin{propo}\label{GIC-Jacobi}
Let $\Lambda_{\alpha, \beta}$ be the commuting 
differential operators defined in Eq.\,$(\ref{MJP-Lambda})$ 
$($with $x$ replaced by $z$$)$ 
related with the Jacobi orthogonal polynomials. 
Assume further that there exists $1\le i\le n$ such that 
$\alpha_i=0$ or $\beta_i=0$. 
Then $\im'\Lambda_{\alpha, \beta}=\cz$, 
and hence Conjecture \ref{GIC}  
holds for $\Lambda_{\alpha, \beta}$.
\end{propo}

\pf Without losing any generality, we may assume that 
$\alpha_1=0$ or $\beta_1=0$. Here we only prove 
the case that $\beta_1=0$. The proof of the case that  
$\alpha_1=0$ is similar.

Under the assumption above, the first component $\Lambda_1$ 
of $\Lambda_{\alpha, \beta}$ in Eq.\,$(\ref{MJP-Lambda})$ is 
the differential operator 
$\Lambda_1= \p_1- \alpha_1(1-z_1)^{-1}$.
Now we apply the change of variables $z_1\to z_1+1$ 
and $z_i\to z_i$ for any $2\le i\le n$, then 
$\Lambda_1$ becomes the differential operator
$\Phi_{\alpha_1}=\p_1+\alpha_1 z_1^{-1}$.

Let $\cA\!:=\bC[z_2, z_3, ..., z_n]$ and view $\Phi_{\alpha_1}$ 
as a differential operator of $\cA[z_1^{-1}, z_1]$. Since $\alpha_1>-1$, 
by applying Corollary \ref{C4.2.3} to $\Phi_{\alpha_1}$, we have 
\begin{align*}
\cA[z_1] \cap \Phi_{\alpha_1}(\cA[z_1])=\cA[z_1].
\end{align*}
Since $\cA[z_1]=\cz$, we have
\begin{align*}
\im'\Phi_{\alpha_1}=\cz \cap \Phi_{\alpha_1}(\cz)=\cz.
\end{align*}

Hence we also have $\im'\Lambda_1=\cz$. 
Since $\im'\Lambda_1\subset \im'\Lambda_{\alpha, \beta}$,
we have $\im'\Lambda_{\alpha, \beta}=\cz$.
\epfv

Next let us consider the differential operators 
$\Lambda_\alpha$ $(\alpha\in (\bR^{>-1})^{\times n}$) 
related with the multi-variable Laguerre orthogonal polynomials. 
Note that, by Eq.\,(\ref{LgrEqs}) and a similar construction for the multi-variable Jacobi polynomials in Example \ref{MJP},
we know that $\Lambda_\alpha$ is given by 
\begin{align}
\Lambda_\alpha=(\Lambda_{\alpha_1}, 
\Lambda_{\alpha_2}, ..., \Lambda_{\alpha_n}),
\end{align}
where, for any $1\le i\le n$, 
\begin{align}\label{Lambda_i}
\Lambda_{\alpha_i}=\p_i +\alpha_iz_i^{-1}-1.
\end{align}
 
\begin{propo}\label{GIC-Laguerre}
Let $\Lambda_\alpha$ $(\alpha\in (\bR^{>-1})^{\times n})$ 
be the commuting differential operators defined above.
Then, $\im'\Lambda_\alpha=\cz$ iff $\alpha_i=0$ 
for some $1\le i\le n$. 
Hence, Conjectures \ref{GIC} holds for $\Lambda_\alpha$ 
under this condition. 
\end{propo}

\pf $(\Leftarrow)$ Without losing any generality, 
we may assume that $\alpha_1=0$.  
Then, by Eq.\,(\ref{LgrEqs}), we have 
$\Lambda_{\alpha_1=0}=\p_1-1$.
Since $\Lambda_{\alpha_1=0}(-1)=1$, we have 
$1\in \im'(\Lambda_{\alpha_1=0})\subset 
\im'\Lambda_\alpha$. Then, 
by Lemma \ref{L4.1.5}, $(c)$, we have 
$\im'\Lambda_\alpha=\cz$.

This result can also be proved by the following 
more straightforward argument  
(without using Lemma \ref{L4.1.5}). 
Note that $\Lambda_{\alpha_1=0}=\frac d{dz_1}-1$ 
is invertible as a linear operator of $\cz$. Its 
inverse operator is given by 
\begin{align*}
\Lambda_{\alpha_1=0}^{-1}=\left(\p_1-1\right)^{-1}=
-1- \sum_{k=1}^{+\infty} \p_1^k.
\end{align*}
Note that the infinity sum on the right hand side 
of the equation above is a well-defined linear map 
of $\cz$. 

Since $\Lambda_{\alpha_1=0}$ is invertible, we have 
$\im' \Lambda_{\alpha_1=0}=\cz$. Hence we also have 
$\im' \Lambda_{\alpha}=\cz$.

$(\Rightarrow)$ Assume that $\im' \Lambda_{\alpha}=\cz$ 
but $\alpha_i\ne 0$ for any $1\le i\le n$.
In particular, we have $1\in \im' \Lambda_{\alpha}$. 
So there exist $h_i(z)\in \cz$ $(1\le i\le n)$
such that 
\begin{align}\label{GIC-Laguerre-pe2}
1=\sum_{i=1}^n \Lambda_{\alpha_i} h_i(z)
=\sum_{i=1}^n  (\p_i +\alpha_i z_i^{-1}-1) h_i(z).
\end{align}

Now we view the RHS of Eq.\,(\ref{GIC-Laguerre-pe2}) above 
as a Laurent polynomial in $z_1$ with coefficients in 
$\bC[z_2^\pm, ..., z_n^\pm]$. Then the coefficient of 
$z_1^{-1}$ of the RHS is given by 
$\alpha_1 h_1(z)|_{z_1=0}$ which, 
by Eq.\,(\ref{GIC-Laguerre-pe2}), must be 
zero. Hence we also have 
$h_1(z)|_{z_1=0}=0$ since 
$\alpha_1\ne 0$.  Therefore, 
$h_1(z)=z_1\tilde h_1(z)$ for some $\tilde h_1(z)\in \cz$. 
Apply similar arguments to $h_i(z)$ for $2\le i\le n$, 
we see that, for any $1\le i\le n$, 
there exists $\tilde h_i(z)\in \cz$ 
such that 
\begin{align}\label{GIC-Laguerre-pe5}
h_i(z)=z_i \tilde h_i(z).
\end{align}

Next, for any fixed $1\le i\le n$ and $f(z)\in \cz$, 
it is easy to check that 
\begin{align}\label{GIC-Laguerre-pe6} 
\p_i(z^\alpha f(z) e^{-\sum_{i=1}^n z_i})= 
z^\alpha (\Lambda_{\alpha_i} f(z) ) e^{-\sum_{i=1}^n z_i}.
\end{align}
Consequently, by Eq.\,(\ref{GIC-Laguerre-pe5}) we have  
\begin{align}\label{GIC-Laguerre-pe7} 
& \int_0^{+\infty} (\Lambda_{\alpha_i} h_i)z^\alpha e^{-\sum_{i=1}^n z_i}\, d z_i 
=\int_0^{+\infty} 
\p_i (z^\alpha h_i e^{-\sum_{i=1}^n z_i})\, d z_i \\
&=(z^\alpha h_i(z)e^{-\sum_{i=1}^n z_i})\left. \right|_{z_i=0}^{z_i=+\infty}
=(z^\alpha z_i \tilde h_i(z)e^{-\sum_{i=1}^n z_i})\left. 
\right|_{z_i=0}^{z_i=+\infty} \nno \\
&=0, \nno 
\end{align}
where the last equality above follows from the fact 
that $\alpha_i+1>0$ since $\alpha_i>-1$.

Combining 
Eqs.\,(\ref{GIC-Laguerre-pe2})--(\ref{GIC-Laguerre-pe7}),  
we have
\begin{align}\label{GIC-Laguerre-pe8} 
 & \int_{ (\bR^{\ge 0})^{\times n} } z^\alpha  e^{-\sum_{i=1}^n z_i} \, d z
=\sum_{i=1}^n \int_{ (\bR^{\ge 0})^{\times n} } 
(\Lambda_{\alpha_i} h_i) z^\alpha  e^{-\sum_{i=1}^n z_i} \, dz  \\
=\sum_{i=1}^n & \int_{(\bR^{\ge 0})^{\times (n-1)} }
\left( 
\int_0^{+\infty} (\Lambda_{\alpha_i} h_i)  z^\alpha e^{-\sum_{i=1}^n z_i}\, d z_i \right) dz_1\cdots \widehat{dz_i}\cdots dz_n \nno \\
=0.\,\,\,\, & \nno 
\end{align}

But this is a contradiction since 
$z^\alpha e^{-\sum_{i=1}^n z_i}$ is continuous and 
positive everywhere over $(\bR^{>0})^{\times n}$.
\epfv

Next, motivated by the differential operators related 
with the multi-variable Hermite polynomials and also 
the differential operators $\Phi_\lambda$ 
$(\lambda\in \bC^n)$ in Subsection \ref{S4.2}, 
we consider the following family of 
commuting differential operators.

For any $\alpha=(\alpha_1, \alpha_2, ..., \alpha_n)\in \bC^n$, we set 
\begin{align}\label{hmt-Psi_i}
\Psi_{\alpha_i}\!:=\p_i +\alpha_i z_i.
\end{align}
and 
\begin{align}\label{hmt-Psi}
\Psi_\alpha\!:=(\Psi_{\alpha_1}, 
\Psi_{\alpha_2}, ..., \Psi_{\alpha_n}).
\end{align}

Note that, when $\alpha_i=-2$ for any 
$1\le i\le n$, the differential 
operators $\Psi_\alpha$ becomes 
the differential operators related 
with the multi-variable Hermite polynomials. 
This can be seen from Eq.\,(\ref{HmtEqs}) 
and a similar construction in Example \ref{MJP} 
for the multi-variable Jacobi polynomials.

Note also that the differential 
operators $\Psi_\alpha$ $(\alpha\in \bC^n)$
are actually the differential operators of 
the polynomial algebra 
$\cz$ itself. So in this case, we have 
$\im'\Psi_\alpha= \im\Psi_\alpha$ for any 
$\alpha\in \bC^n$.

\begin{propo}\label{IC-Hermite}
Let $\Psi_\alpha$ $(\alpha\in \bC^n)$ be defined above.
Then, $\im \Psi_\alpha=\cz$ iff $\alpha_i=0$ 
for some $1\le i\le n$. Hence Conjectures \ref{IC} 
and \ref{GIC} hold under this condition.
\end{propo}

\pf $(\Leftarrow)$ Assume that $\alpha_i=0$ 
for some $1\le i\le n$.  
Then, by Eq.\,(\ref{hmt-Psi_i}), we have 
$\Lambda_{\alpha_i=0}=\p_i$ which is obviously 
a surjective linear map from $\cz$ to $\cz$.
Hence we have $\im \Psi_{\alpha_i=0}=\cz$ and  
$\im \Psi_{\alpha}=\cz$.

$(\Rightarrow)$ Assume that 
$\im \Psi_{\alpha}=\cz$ but $\alpha_i\ne 0$ 
for any $1\le i\le n$. We derive a contradiction 
as follows.

First, by applying the change of variables 
$z_i \to \sqrt{-2/\alpha_i}\, z_i$ 
$(1\le i\le n)$, the differential 
operators $\Psi_{\alpha_i}$ in Eq.\,(\ref{hmt-Psi_i}) 
becomes $\sqrt{-\alpha_i/2}\, (\p_i-2z_i)$. 
Set $\Lambda_i\!:=\p_i-2z_i$ and 
$\Lambda=(\Lambda_1, \Lambda_2,..., \Lambda_n)$.
From the argument above, it is easy see that 
$\im\Lambda=\im\Psi_\alpha=\cz$.

Let $h_i(z)\in \cz$ $(1\le i\le n)$
such that 
\begin{align}\label{IC-Hermite-pe1}
1=\sum_{i=1}^n \Lambda_i h_i.
\end{align}

Note that, for any fixed $1\le i\le n$,  
it is easy to check that 
\begin{align}\label{IC-Hermite-pe2} 
\p_i( h_i e^{-\sum_{i=1}^n z_i^2})= 
(\Lambda_i h_i ) e^{-\sum_{i=1}^n z_i^2}.
\end{align}
Consequently, we have  
\begin{align}\label{IC-Hermite-pe3}
& \int_\bR (\Lambda_i h_i) e^{-\sum_{i=1}^n z_i^2}\, d z_i 
=\int_\bR  
\p_i (h_i e^{-\sum_{i=1}^n z_i^2})\, d z_i \\
&=\left. ( h_i(z) e^{-\sum_{i=1}^n z_i^2}) \right|_{z_i=-\infty}^{z_i=+\infty}
=0. \nno 
\end{align}

Then, by using Eqs.\,(\ref{IC-Hermite-pe1})--(\ref{IC-Hermite-pe3}) 
and applying a similar argument as in Eq.(\ref{GIC-Laguerre-pe8}), 
we have 
\begin{align*}
\int_{\bR^n} e^{-\sum_{i=1}^n z_i^2}\, dz=\sum_{i=1}^n \int_{\bR^n } 
(\Lambda_i h_i)  e^{-\sum_{i=1}^n z_i^2} \, dz  =0, 
\end{align*}
which is a contradiction since 
$e^{-\sum_{i=1}^n z_i^2}$ is continuous 
and positive everywhere on $\bR^n$.
\epfv 

Now let us consider some cases of Conjecture \ref{GMC}.

\begin{propo}\label{FiniteCase}
Conjecture \ref{GMC} holds for any open $B\subset \bR^n$ with 
any atomic measure $\sigma$ which is supported at 
finitely many points of $B$.
\end{propo}

\pf Let $S=\{u_1, u_2, ..., u_k\}\subset B$ be the 
support of $\sigma$, i.e. $\sigma(u_i)> 0$ 
$(1\le i\le k)$ and, for any measurable subset 
$U\subset B$, we have
\begin{align*}
\sigma(U)=\sum_{u \in S\cap U} \sigma(u).
\end{align*}

Note first that, for any $f(z)\in \cz$, $f(z)\in \cM_B(\sigma)$ 
iff
\begin{align}\label{FiniteCase-pe2}
\int_B f(z)\, d\sigma=\sum_{i=1}^k f(u_i) 
\sigma(u_i)=0.
\end{align}

Therefore, for any $f(z)\in \cz$ with  
$f^m(z)\in \cM_B(\sigma)$ for any $m\ge 1$, we have
\begin{align}\label{FiniteCase-pe3}
\int_B f^m(z)\, d\sigma=\sum_{i=1}^k f^m(u_i) 
\sigma(u_i)=0
\end{align}
for any $m\ge 1$

If $f(u_i)=0$ for all $1\le i\le k$, then, 
for any $g(z)\in \cz$ and $m\ge 1$, we also have 
$(f^m g)(u_i)=0$ for any $1\le i\le k$.  
Hence, for any $m\ge 1$, $f^mg$ also satisfies 
Eq.\,(\ref{FiniteCase-pe2}) 
and lies in $\cM_B(\sigma)$. Therefore, 
Conjecture \ref{GMC} holds for $f(z)$.

Assume $f(u_i)$ $(1\le i\le k)$ are not all zero. 
Let $\{c_1, c_2, ..., c_s\}$ be the set of
all distinct nonzero values of $f(z)$ attained 
over $S$. For any $1\le j\le s$, let $S_j$ be the 
subset of elements of $u\in S$ such that 
$f(u)=c_j$ and $k_j$ the cardinal number 
of $S_j$. Then, Eq.\,(\ref{FiniteCase-pe3}) 
can be re-written as 
\begin{align}\label{FiniteCase-pe4}
0=\sum_{i=1}^k f^m(u_i) 
\sigma(u_i)=\sum_{j=1}^s c_j^m 
\sum_{u\in S_j} \sigma(u),
\end{align}
Since the equation above holds 
for any $m\ge 1$, by using the invertibility 
of the Vandermonde matrices, 
it is easy to check that  
\begin{align}\label{FiniteCase-pe5}
0=\sum_{u\in S_j} \sigma(u_j).
\end{align}
for any $1\le j\le s$.

But this is a contradiction since $\sigma(u_i)>0$ for any 
$1\le i\le k$. \epfv 

Next, let us consider Conjecture \ref{GMC} for 
the Jacobi orthogonal polynomials. 
Contrast to Conjecture \ref{GIC} for the Jacobi polynomials 
(See Proposition \ref{GIC-Jacobi}), the only case 
of Conjecture \ref{GMC} that we know is   
the case of the one-variable Jacobi polynomials 
with $\alpha=\beta=0$. In this case the weight function 
$w(z)\equiv 1$. Note that, by Remark \ref{R4.1.3}, $(b)$, 
the Jacobi polynomials in this case are actually the 
Legendre polynomials. 

\begin{propo}  
Let $a, b\in \bR$ with $a>b$ and $d\sigma=dz$. Then 
Conjecture \ref{GMC} holds for open interval 
$B\!:=(a, b)\subset \bR$ with the Lebesgue measure 
$\sigma$. 
\end{propo}

The proposition above follows directly from the 
following theorem.

\begin{theo}\label{Folk}
Let $a<b\in \bR$ and $f(z)\in \bC[z]$. Assume that, there exists 
$N>0$ such that $\int_a^b f^m(z)\, dz=0$ for any $m\ge N$. 
Then $f(z)=0$.
\end{theo}

It seems that the theorem above is known but we failed to find 
any published proof in the literature. We did notice that 
Madhav V. Nori \cite{N} has studied the problem above in a much more general setting. It is very possible that Theorem \ref{Folk} will follow from some results obtained in \cite{N}.

Jean-Philippe Furter \cite{FZ} informed the author that he and his colleague Changgui Zhang have got an analytic proof for Theorem \ref{Folk}, which is under preparation. Mitya Boyarchenko \cite{B} also sent the author a sketch of his brilliant but unpublished proof. Surprisingly, Boyarchenko's proof is purely algebraic and uses only some results from algebraic number theory such as Dirichlet's theorem on arithmetic progressions, etc. 

Next, we end this section with a connection of the one-variable case of Conjecture \ref{GMC} with the so-called {\it polynomial moment problem} proposed by M. Briskin, J.-P. Francoise and Y. Yomdin in the series of papers \cite{BFY1}-\cite{BFY5}. 
The {\it polynomial moment problem} is mainly motivated by the center problem for the complex Abel equation. The problem was  recently solved by F. Pakovich and M. Muzychuk \cite{PM} 
(See the theorem below). For more details on the 
{\it polynomial moment problem}, see the references quoted above and also citations therein. The author is very grateful to Harm Derksen, Jean-Philippe Furter, Jeffrey C. Lagarias, Leonid Makar-Limanov, Lucy Moser-Jauslin for communications and suggestions on this connection, and also to Fedor Pakovich for communications on his joint work \cite{PM} with Mikhail Muzychuk. 
 
Recall that the {\it polynomial moment problem} is the following problem: given any polynomial $f(z)\in \cz$ and $a\ne b\in \bC$,  
find all polynomials $q(z)\in \cz$ such that, for any $m\ge 0$, 
\begin{align}\label{PMP-e1}
\int_a^b f^m(z) q(z)dz=0.
\end{align}

The problem above was solved recently by the following theorem 
obtained by F. Pakovich and M. Muzychuk \cite{PM}.
 
\begin{theo} $(${\bf Pakovich and Muzychuk}$)$ 
Let $a\neq b\in \bC$ and $f(z)\in \cz$. A non-zero polynomial 
$q(z)\in \cz$ satisfies Eq.\,$($\ref{PMP-e1}$)$ for any $m\ge 0$ 
iff there exist some polynomials $Q_j(z)$, $f_j (z)$ 
and $W_j(z)$ $(j\in J)$ such that 

$(1)$ $W_j(a)=W_j(b)$ for any $j\in J$;

$(2)$ $q(z)=\sum_{j\in J} Q_j'(W_j(z))\, W_j'(z)$;

$(3)$ $f(z) = f_j(W_j (z))$ for any $j\in J$.
\end{theo}

Note that, with the same notation as in Pakovich and   Muzychuk's theorem above, if we choose $a<b$, $B=(a, b)\subset \mathbb R$ and the measure $d\sigma=q(z)dz$ (ignoring the positivity requirement  on the measure $\sigma$ for a moment), then Pakovich and Muzychuk's theorem above gives all polynomials $f(z)\in \cz$ 
such that $\int_B f^m(z)\, d\sigma=0$ for any $m\ge 0$.

But, unfortunately, Pakovich and Muzychuk's theorem requires 
the integral in Eq.\,(\ref{PMP-e1}) vanish when $m=0$, 
i.e. $\int_B d\sigma = \int_B q(z) dz=0$. From this requirement, it is easy to see that $d\sigma=q(z)dz$ can not be a positive measure 
on the interval $B=(a, b)$. 
For example, $q(z)$ can not be any nonzero polynomial with real coefficients such that $q(c)\ge 0$ for any $c\in B$. 
Therefore, Pakovich and Muzychuk's theorem can not be applied directly to approach Conjecture \ref{GMC}. 

Nevertheless, it is still very interesting to see if some of the  techniques (instead of the main theorem) in \cite{PM} can somehow be applied to study Conjecture \ref{GMC} for the cases when $q(z)$ are polynomials or analytic  functions in one variable which are non-negative over some open intervals 
of the real line.

On the other hand, we see that Conjecture \ref{GMC} also raises a new question on the {\it polynomial moment problem}, namely, what is the solution of the {\it polynomial moment problem} with the (slightly weaker) condition that the integrals in Eq.\,(\ref{PMP-e1}) vanish for any $m\ge 1$ but not necessarily for $m=0$? We believe this question is also very interesting 
to investigate.

\renewcommand{\theequation}{\thesection.\arabic{equation}}
\renewcommand{\therema}{\thesection.\arabic{rema}}
\setcounter{equation}{0}
\setcounter{rema}{0}

\section{\bf Some General Results on Mathieu Subspaces}
\label{S5}

Note that the Mathieu conjecture (Conjecture \ref{MC}), 
the \IC (Conjecture \ref{IC}) and its generalizations 
(Conjectures \ref{GIC} and \ref{GMC})
discussed in the previous sections are all  
about whether or not certain subspaces are 
Mathieu subspaces of their ambient commutative algebras. 
Furthermore, this is also the case for the well-known Jacobian conjecture and the Dixmier conjecture through their equivalences (See \cite{IC} for the discussions on these equivalences) to some special cases of the {\bf IC}. Therefore, it is necessary and important to study Mathieu subspaces separately 
in a more general setting.

In this section, we give some examples of Mathieu subspaces from other sources and derive some general results on this newly introduced notion.

First, let us generalize the notion of Mathieu subspaces defined in 
Definition \ref{M-ideal} to associative but not necessarily 
commutative algebras.

\begin{defi}
Let $\cA$ be an associative algebra over a commutative 
ring $R$ and $\cM$ a $R$-subspace of $\cA$. 
We say that $\cM$ is a {\it left} Mathieu subspace of $\cA$ 
if the following property holds: 
for any $a, b\in \cA$ with $a^m\in \cM$ for any 
$m\ge 1$, there exists $N\ge 1$ $($depending on $a$ and $b$$)$ such that $b a^m \in \cM$ for any $m\ge N$.
\end{defi}

We define {\it right} Mathieu subspaces and 
also ({\it two-sided}) Mathieu subspaces in 
the obvious ways.

It is easy to see that any left ideal is automatically a left Mathieu subspace. Similarly, this is also the case for right and two-sided ideals. Therefore, the notion of Mathieu subspaces can be viewed as a generalization of the notion of ideals even for noncommutative algebras. But, as we will see from examples to be  discussed in this section, many Mathieu subspaces are not ideals. Actually they are not even closed under the product 
of the ambient algebras.

We start with the following noncommutative examples of Mathieu subspaces. 

\begin{exam}
For any $n\ge 1$ and integral domain $R$ of characteristic zero, let $\cA=M_{n\times n}(R)$ be the algebra of $n\times n$ matrices with entries in $R$ and $\cM$ the subspace of trace-zero matrices.
 
Note that, for any $A \in \cA$, $A^m\in \cM$ for any $m\ge 1$ 
iff $A$ is nilpotent. Then, for any $B\in \cA$, we have
$BA^m=A^mB=0\in \cM$ for any $m\ge n$.  Therefore, $\cM$ is 
a two-sided Mathieu subspace of $\cA$ but certainly can not be an ideal of $\cA$ unless $n=1$.
\end{exam}

Next, from additive valuations on polynomial algebras, 
we can get the following family of Mathieu subspaces 
of polynomial algebras. 

\begin{exam}
For any $n\ge 1$ and any integral domain $R$, let 
$\cA=R[z]$ be the polynomial algebra over 
$R$ in $n$ variables $z$. 
For any linear functional  $\nu: \bR^n \to \bR$, we define an additive valuation $\mbox{ord}_\nu: \cA\to \bR\cup\{+\infty\}$ by setting, $\mbox{ord}_\nu (0)\!:=+\infty$ 
and, for any $0\neq f(z) \in \cA$,
\begin{align} 
\mbox{ord}_\nu (f):=\min \{\nu(\alpha) \,|\, \mbox
{the coefficient of $z^\alpha$ in $f(z)$ is not zero}\}. 
\end{align} 
For any $c\in \bR$, let $\cM_c$ be the subspace 
of polynomials $f\in \cA$ such that 
$\mbox{ord}_\nu(f) \ge c$.

Then it is easy to check that, for any $c>0$, 
$\cM_c$ is a Mathieu subspace of $\cA$ 
but not necessarily an ideal of $\cA$ 
if $\mbox{ord}_\nu (f)<0$ 
for some $f\in \cA$. 
\end{exam}

More generally, we have the following family of 
examples from commutative rings 
(viewed as $\bZ$-algebras) 
with valuations.

\begin{exam}
Let $\cA$ be any commutative ring and $\nu$ a real-valued 
$($additive$)$ valuation 
$($\cite{AM}, \cite{Msu}$)$ of $\cA$, i.e. 
$\mu: \cA\to \bR\cup\{+\infty\}$ such that, for any $x, y \in \cA$, we have 
\begin{align}
\nu(x)&=+\infty \mbox{ iff } x=0,\\
\nu(xy) &=\nu(x)+ \nu(y),\\
\nu(x+y)&\ge \min\{\nu(x), \nu(y)\}.
\end{align}

For any $c\in \bR$, let $\cM_c$ be the subspace of elements 
$x$ of $\cA$ with $\nu(x)\ge c$.
Then it is easy to check that, 
for any $c>0$, $\cM_c$ is a 
Mathieu subspace of $\cA$ but not necessarily 
an ideal of $\cA$. 
\end{exam}

Similar as for ideals, we say a Mathieu subspace $\cM$ of 
an algebra $\cA$ is {\it proper} if $\cM\ne 0$ and $\cM\ne \cA$. 
From the example above, we see that some fields may actually have some proper Mathieu subspaces. Also, unlike proper ideals, proper Mathieu subspaces may contain some invertible elements. 
But, as we will see below, proper left or right 
Mathieu subspaces can not contain
the identity element of the 
ambient algebras.  

\begin{lemma}\label{No-1-lemma}
Let $\cA$ be any algebra and $\cM$ a proper left or right Mathieu subspace of $\cA$. Then, $1\not \in \cM$. 
\end{lemma}
\pf Assume otherwise. Then, for any $m\ge 1$, 
$1^m=1\in \cM$. Then, for any $b\in \cA$, 
$b=1^mb=b1^m \in \cM$ when $m>>0$. 
Hence, we have $\cM=\cA$ which 
is a contradiction.
\epfv 

The next proposition will give us more examples 
of Mathieu subspaces of polynomial algebras, 
which are not necessarily ideals.

\begin{propo}\label{FiniteCase2}
Let $K$ be a field of any characteristic.  
For any finite subset $S=\{u_1, u_2, ..., u_k\}\subset K^n$ 
and $\sigma=\{a_1, a_2, ..., a_k\} \subset K^{\times}$, 
denote by $\cM(\sigma)$ the subspace of polynomials 
$f(z)\in K[z]$ such that
\begin{align}\label{FiniteCase2-e1}
\sum_{i=1}^k a_i f(u_i)=0.
\end{align}

Then, $\cM(\sigma)$ is a 
Mathieu subspace of $K[z]$ iff, 
for any non-empty subset $S'\subset S$, 
\begin{align}\label{FiniteCase2-e2}
\sum_{u\in S'} u \ne 0.
\end{align}
\end{propo}

\pf First, the $(\Leftarrow)$ part of the proposition 
can be proved by similar arguments as in the proof 
of Proposition \ref{FiniteCase}. So we skip it here.

To show the $(\Rightarrow)$ part, assume that 
there exists a non-empty $S'\subset S$ such that 
Eq.\,(\ref{FiniteCase2-e2}) holds. 

Let $f(z)\in \Kz$ such that
\begin{align}\label{FiniteCase2-pe1}
f(u_i)=
\begin{cases}
1 &\mbox{ if } u_i\in S'\\
0 &\mbox{ if } u_i\not \in S'. 
\end{cases}
\end{align}

Note that, such a polynomial $f(z)$ always exists. For example, when $n=1$ (the idea for $n>1$ is similar), we may choose 
\begin{align*} 
f(z)= \sum_{u\in S'} \frac {\prod_{c\in S; c\ne u} (z-c)} 
{\prod_{c\in S; c\ne u} (u-c)} 
\end{align*}

For any $m\ge 1$, we have
\begin{align*}
\sum_{i=1}^k a_i f^m(u_i)=\sum_{i\in S'} a_i = 0.
\end{align*}
Therefore, we have that 
$f^m(z)\in \cM(\sigma)$ for any $m\ge 1$.

Now fix an $u_j \in S'$ and
$g(z)\in K[z]$ such that $g(u_j)=1$ and 
$g(u_i)=0$ for any $1\le i\ne j\le k$. 
Then, for any $m\ge 1$, 
we have
\begin{align*}
\sum_{i=1}^k a_i (f^m g)(u_i)=\sum_{i=1}^k a_i f^m(u_i) 
g(u_i)=a_j g(u_j)=a_j\ne 0.
\end{align*}

Therefore, $f^mg \not \in\cM(\sigma)$ for 
any $m\ge 1$. Hence $\cM(\sigma)$ is not 
a Mathieu subspace of $K[z]$.
\epfv

\begin{rmk}
Note that, when $k=1$, $\cM(\sigma)$ in 
Proposition \ref{FiniteCase} 
is always an ideal of $\cz$. 
But if $k>1$, $\cM(\sigma)$ in general 
is not an ideal. For example, 
for any $k\ge 2$, choose $K=\bC$, $a_i=1$ and 
any distinct $u_i\in \bC$ $(1\le i\le k)$.
\end{rmk}

For Laurent polynomial algebras $\czz$, by Duistermaat and van der Kallen's theorem, Theorem \ref{ThmDK}, we see that the subspace of Laurent polynomials with no constant terms is a Mathieu subspace of $\czz$. Another example of Mathieu subspaces of $\czz$ is given by the following theorem which was first conjectured in \cite{GVC} and later was proved in \cite{EWZ}.

\begin{theo}
Let $\cM$ be the subspace of $\czz$ of Laurent polynomials $f(z)$
with no holomorphic part, i.e. $[z^\alpha]f(z)=0$ for any 
$\alpha\in\bN^n$. Then, $\cM$ is a Mathieu subspace of $\czz$. 
\end{theo}
 
Next we show that some properties of ideals  are also shared by Mathieu subspaces. 

\begin{propo}\label{Intersection}
Let $\cA$ and $\cB$ be any algebras over a commutative ring $R$, and $\varphi:\cA\to \cB$ an $R$-algebra homomorphism. Then,

$(a)$ for any left Mathieu subspaces $\cM_i$ 
$(1\le i\le m)$ of $\cA$, $\cap_{1\le i\le m} \cM_i$ 
is also a left Mathieu subspace of $\cA$.

$(b)$ for any left Mathieu subspace $\cN$ of $\cB$, 
$\varphi^{-1}(\cN)$ is also a left Mathieu 
subspace of $\cA$.
\end{propo}
\pf $(a)$ Set $\cM\!:=\cap_{1\le i\le m} \cM_i$. 
It is easy to see that $\cM$ is also a 
$R$-subspace of $\cA$.

Let $a, b\in \cA$ with $a^m\in \cM$ for 
any $m\ge 1$. For any $1\le i\le m$, let $N_i\in\bN$ such that 
$b a^m\in \cM_i$ for any $m\ge N_i$. Let 
$N\!:=\max\{ N_i\,|\, 1\le i\le m\}$. Then, for
any $m\ge N$, we have $ba^m\in \cM_i$ for  
any $1\le i\le m$. Hence $ba^m\in \cM$  
for any $m\ge N$.

$(b)$ Again, it is easy to see that $\varphi^{-1}(\cN)$ is a 
$R$-subspace of $\cA$ since $\varphi$ is an $R$-algebra 
homomorphism.

Let $a, b\in \cA$ with
$a^m\in \varphi^{-1}(\cN)$ for any 
$m\ge 1$. Set $x\!:=\varphi(a)$ and 
$y\!:=\varphi (b)$. Then, for any $m\ge 1$, we have
$x^m=\varphi^m(a)=\varphi(a^m) \in \cN$. 
Since $\cN$ is a left Mathieu 
subspace of $\cB$, there exists $N\in \bN$ such that 
$yx^m \in \cN$ for any $m\ge N$. But 
$yx^m=\varphi(b)\varphi^m(a) =\varphi(ba^m)$, so
we have $ba^m\in \varphi^{-1}(\cN)$ for
any $m\ge N$. Therefore, $\varphi^{-1}(\cN)$ is a left Mathieu subspace of $\cA$.
\epfv

\begin{rmk}
From the proof above, it is easy to see that Proposition 
\ref{Intersection} also holds for right or two-sided 
Mathieu subspaces. 
\end{rmk}

\begin{propo}\label{M-extention}
Let $K$ be any field with uncountably many elements and 
$\cA$ a commutative $K$-algebra. Let $z=(z_1, z_2, ..., z_n)$ 
be $n$ free commutative variables. Then, for any Mathieu subspace 
$\cM$ of $\cA$, $\cM[z]\subset \cAz$ is a Mathieu subspace of 
$\cAz$.
\end{propo}

\pf Note first that, by using induction on the number of free 
variables $z$, it will be enough to show the proposition for 
the one-variable case. So we assume $n=1$.

We use the contradiction method. Assume the proposition is false. Then there exist $f(z), g(z)\in K[z]$ with $f^m\in \cM[z]$ 
for any $m\ge 1$ but $f^k(z)g(z)\not \in \cM[z]$ for infinitely 
many positive integers $k$. Let $\{ m_i\in \bN\}$ be a strictly increasing sequence of positive integers such that 
$f^{m_i}g \not \in \cM[z]$ for any $i\ge 1$. 

Note that, for any $h(z)\in \cAz$ of degree 
$d\!:=\deg h(z)\ge 0$ and $h(z)\not \in \cM[z]$, 
by using the invertibility of 
the Vandermonde matrices, it is easy to check that 
there are at most $d$ distinct elements $c\in K^\times$ 
such that $h(c) \in \cM$. Otherwise, all the coefficients 
of $h(z)$ would be in $\cM$ and $h(z)\in \cM[z]$.

Therefore, for each fixed $m_i$, 
there are only finitely many 
$c\in K$ such that $f^{m_i}(c)g(c)\in \cM$.
Since $K$ has uncountably many distinct elements, 
there exists $b\in K$ such that 
$f^{m_i}(b)g(b)\not \in \cM$ 
for all $i\ge 1$.

But, on the other hand, since $f^{m}(z) \in \cM[z]$ 
for any $m\ge 1$, all the coefficients of $f^{m}(z)$ 
are in $\cM$. Since $\cM$ is a $K$-subspace of $\cA$, 
we have $f^{m}(b)\in \cM$ for any $m\ge 1$. 
Furthermore, since $\cM$ is a Mathieu subspace of $\cA$ and 
$f(b)^m=f^m(b) \in \cM$ for any $m\ge 1$, 
we have $f^m(b)g(b)\in \cM$ when $m>>0$. 
In particular, $f^{m_i}(b)g(b)\in \cM$ when $i>>0$. 
Hence we get a contradiction.
\epfv

Finally, one remark on the sums of Mathieu subspaces is as follows.

Note that, by Proposition \ref{Intersection}, $(a)$, the intersection of any finitely many Mathieu subspaces 
is always a Mathieu subspace. Naturally, one may wonder 
if the sum of any finitely many Mathieu subspaces is also a Mathieu subspace. But this is not true in general. 

\begin{exam}
Let $\cA=\cz$ in one variable $z$. Let $\cM_1$ and $\cM_2$ be the one-dimensional subspace spaces of $\cz$ spanned by $1+z$ and $1-z$, respectively. Then, it is easy to check that both 
$\cM_1$ and $\cM_2$ are Mathieu subspaces of $\cz$ 
and $\cM\!:=\cM_1+\cM_2=\bC\cdot 1 +\bC \cdot z$. But, 
since $1\in \cM$ and $\cM\ne \cz$, by Lemma \ref{No-1-lemma}  
$\cM$ is not a Mathieu subspace of $\cz$.
\end{exam}

{\small \sc Department of Mathematics, Illinois State University,
Normal, IL 61790-4520.} 

{\em E-mail}: wzhao@ilstu.edu.


\begin{thebibliography}{FLM2}


\bibitem[AS]{AS} {\it Handbook of mathematical functions with formulas, graphs, and mathematical tables.} Edited by Milton Abramowitz and Irene A. Stegun. Reprint of the 1972 edition. 
Dover Publications, Inc., New York, 1992. [MR0757537].

\bibitem[AE]{AE} P. K. 
Adjamagbo and A. van den Essen, {\it 
A proof of the equivalence of the Dixmier, Jacobian and Poisson conjectures.} Acta Math. Vietnam. 32 (2007), no. 2-3, 205--214. [MR2368008]. 

\bibitem[AM]{AM} M. F. Atiyah and I. G. Macdonald, {\it Introduction to commutative algebra}.
Addison-Wesley Publishing Co., 1969. [MR0242802].

\bibitem[BCW]{BCW} H. Bass, E. Connell, D. Wright, {\it The Jacobian
conjecture, reduction of degree and formal expansion of the inverse}. Bull.  Amer. Math.  Soc.  \textbf{7}, (1982), 287--330. 
[MR 83k:14028]. Zbl.539.13012.


\bibitem[BK]{BK} A. Belov-Kanel and M. Kontsevich, {\it 
The Jacobian conjecture is stably equivalent to the Dixmier conjecture.} (English, Russian summary) 
Mosc. Math. J. {\bf 7} (2007), no. 2, 209--218, 349. 
[MR2337879].  

\bibitem[B]{B} M. Boyarchenko, {\it Personal communications}.

\bibitem[BFY1]{BFY1} M. Briskin, J.-P. Francoise, Y. Yomdin, 
{\it Une approche au probleme du centre-foyer de
Poincare}. C. R. Acad. Sci., Paris, Ser. I, Math. {\bf 326} (1998), No.11, 1295-1298. [MR1649140].

\bibitem[BFY2]{BFY2} M. Briskin, J.-P. Francoise, Y. Yomdin, {\it Center conditions, compositions of polynomials and
moments on algebraic curve}. Ergodic Theory Dyn. Syst. {\bf 19} (1999), no 5, 1201--1220. [MR1721616].

\bibitem[BFY3]{BFY3} M. Briskin, J.-P. Francoise, Y. Yomdin, {\it Center condition II: Parametric and model center
problems}. Isr. J. Math. {\bf 118} (2000), 61--82. [MR1776076].

\bibitem[BFY4]{BFY4} M. Briskin, J.-P. Francoise, Y. Yomdin, {\it Center condition III: Parametric and model center
problems}. Isr. J. Math. {\bf 118} (2000), 83--108. [MR1776077].

\bibitem[BFY5]{BFY5} M. Briskin, J.-P. Francoise, Y. Yomdin, {\it Generalized moments, center-focus conditions and
compositions of polynomials}. Operator theory, system theory and related topics (Beer-Sheva/Rehovot, 1997), 161--185, Oper. Theory Adv. Appl., {\bf 123} (2001). [MR1821911].

\bibitem[C]{C} T. S. Chihara, {\it An introduction to orthogonal polynomials.} Mathematics and its Applications, Vol. 13. 
Gordon and Breach Science Publishers, 
New York-London-Paris, 1978. [MR0481884].

\bibitem[D]{D} J. Dixmier, {\it Sur les alg\`ebres de Weyl.} Bull. Soc. Math. France, 96 (1968), 209--242. [MR0242897]. 

\bibitem[DK]{DK} J. J. Duistermaat and W. van der Kallen, 
{\it Constant terms in powers of a Laurent polynomial.} 
Indag. Math. (N.S.) {\bf 9} (1998), no. 2, 221--231. [MR1691479].

\bibitem[DX]{DX} C. Dunkl and Y. Xu, {\it Orthogonal polynomials of several variables.} Encyclopedia of Mathematics and its Applications, 81. Cambridge University Press, Cambridge, 2001. [MR1827871].

\bibitem[E]{E} A. van den Essen, {\emph Polynomial automorphisms and the Jacobian conjecture}.  Progress in Mathematics, 190. Birkh\"auser Verlag, Basel, 2000. [MR1790619]. 

\bibitem[EWZ]{EWZ} A. van den Essen, R. Willems and W. Zhao, 
{\it Some results on the vanishing conjecture of differential 
operators with constant coefficients}. 
arXiv:0903.1478[math.AC]. 
Submitted.

\bibitem[FZ]{FZ} J.-P. Furter and C. Zhang, {\it Personal communications}. 

\bibitem[Ke]{Ke} O. H. Keller, 
{\it Ganze Gremona-Transformationen}, 
Monats. Math. Physik {\bf 47} (1939), no.\,1, 299-306. [MR1550818].

\bibitem[Ko]{Ko} T. H. Koornwinder, {\it Two-variable analogues of the classical orthogonal polynomials.} Theory and application of special functions (Proc. Advanced Sem., Math. Res. Center, Univ. Wisconsin, Madison, Wis., 1975), pp.\,435--495. Math. Res. Center, Univ. Wisconsin, Publ. No. 35, Academic Press, New York, 1975. [MR0402146].

\bibitem[Ma]{Ma} O. Mathieu, {\it Some conjectures about invariant theory and their applications.} Alg\`ebre non commutative, groupes quantiques et invariants (Reims, 1995), 263--279, S\'emin. Congr., 2, Soc. Math. France, Paris, 1997. [MR1601155].

\bibitem[Mat]{Msu} H. Matsumura,  {\it Commutative ring theory}.
 2nd edition.  Cambridge University Press, 1989. [MR1011461].

\bibitem[N]{N} M. V. Nori, {\it The integral of powers of a function}. Symposium in Honor of C. H. Clemens (Salt Lake City, UT, 2000), 163--175, {\it Contemp. Math.}, {\bf 312} (2002). 
[MR1941581]. 


\bibitem[PM]{PM}  F. Pakovich and M. Muzychuk, 
{\it Solution of the polynomial moment problem.} To appear in {\it Proc. Lond. Math. Soc.}. See also arXiv:0710.4085v2 [math.CV].

\bibitem[SHW]{SHW} J. A. Shohat; E. Hille and J. L. Walsh, 
{\it A Bibliography on Orthogonal Polynomials.} 
Bull. Nat. Research Council, no. 103. National Research Council of the National Academy of Sciences, Washington, D. C., 1940. [MR0003302].

\bibitem[Si]{Si1} B. Simon, {\it Orthogonal polynomials on the unit circle. Part 1. Classical theory.} American Mathematical Society Colloquium Publications, 54, Part 1. American Mathematical Society, Providence, RI, 2005. [MR2105088].

\bibitem[Sz]{Sz} G. Szeg\"o, {\it Orthogonal Polynomials.} 4th edition. American Mathematical Society, Colloquium Publications, Vol. XXIII. American Mathematical Society, Providence, R.I., 1975. [MR0372517].

\bibitem[T]{Ts} 
Y. Tsuchimoto, {\it Endomorphisms of Weyl algebra and 
$p$-curvatures.} Osaka J. Math. 42 (2005), no. 2, 435--452. 
[MR2147727]. 

\bibitem[Z1]{HNP} W. Zhao, {\it Hessian Nilpotent Polynomials and the Jacobian Conjecture}, Trans. Amer. Math. Soc. 359 (2007), no. 1, 249--274 (electronic).
[MR2247890]. See also math.CV/0409534.

\bibitem[Z2]{GVC} W. Zhao, {\it A Vanishing Conjecture on Differential Operators with Constant Coefficients}, Acta Mathematica Vietnamica, vol 32 (2007), no. 3, 107--134. [MR2368014]. See also arXiv:0704.1691v2 [math.CV].

\bibitem[Z3]{IC} W. Zhao, {\it Images of commuting  Differential Operators of Order One with Constant Leading Coefficients}.  arXiv:0902.0210 [math.CV]. Submitted. 

\end{thebibliography}
\end{document}